\documentclass[leqno, a4paper, 11pt]{article}

\usepackage{amsfonts, amssymb, calrsfs}
\usepackage[dvips]{graphicx}
\usepackage[dvips]{color}
\usepackage{epsfig}
\usepackage{latexsym}

\usepackage[german]{babel}

\def\Z{\mathbf{Z}}

\def\F{\mathbf{F}}
\def\G{\mathbf{G}}
\def\P{\mathbf{P}}

\def\O{\mathcal{O}}
\def\m{\mathfrak{m}}

\def\X{{\mathcal{X}}}
\def\til#1{\widetilde{#1}}
\def\Spec{\mathop{\mathrm{Spec}}\nolimits}

\def\Proj{\mathop{\mathrm{Proj}}\nolimits}

\def\Aut{\mathop{\mathrm{Aut}}\nolimits}
\def\Hom{\mathop{\mathrm{Hom}}\nolimits}

\def\ker{\mathop{\mathrm{ker}}\nolimits}

\def\Gal{\mathop{\mathrm{Gal}}\nolimits}

\def\GL{\mathop{\mathrm{GL}}\nolimits}

\def\longhookrightarrow{\lhook\joinrel\longrightarrow}

\def\bigdownarrow{\vphantom{\bigg|}\Big\downarrow}

\def\pf{{\indent\textit{Beweis.}\ }}
\def\qed{\hfill$\square$}
\newcounter{para}[section]
\renewcommand{\thepara}{\thesection.\arabic{para}}
\renewcommand{\thesection}{\arabic{section}}
\renewcommand{\paragraph}{\refstepcounter{para}
\indent{\bf{\thepara}}}
\def\sectioning#1{\vv 

 \refstepcounter{section}
\indent {\bf \thesection. #1}}
\newcommand{\sd}{\rtimes}
\newcommand{\vv}{\vspace{1ex}}

\newcommand{\hx}{\hat{x}}
\newcommand{\hy}{\hat{y}}
\newcommand{\hz}{\hat{z}}

\umlautlow

\begin{document}

{\footnotesize \noindent Kolumnentitel: Entartung von Gruppenoperationen auf Kurven \hfill 30.\ August 2003 

\noindent Math.\ Subj.\ Class.\ (2000): 14H10, 14H30, 14H37}

\vspace{4cm}

\begin{center} 

{\Large Zur Entartung gez\"ugelter Gruppenoperationen auf Kurven}

\vv

{\sl von} Gunther Cornelissen {\sl und} Fumiharu Kato

\vv

\end{center}

\vv

\vv

\noindent {\small {\bf Zusammenfassung.} 
Man nennt eine Gruppenoperation einer endlichen Gruppe auf einer glatten projektiven Kurve \"uber einen algebraisch abgeschlossenen K\"orper der positiven Charakteristik {\sl gez\"ugelt}, falls alle {\sl zweiten} Verzweigungsgruppen trivial sind (es ist z.B.\ jede Gruppenoperation auf einer ordin\"aren Kurve gez\"ugelt). 
Wenn die Verzweigungsindizes gewisse numerische Bedingungen erf\"ullen, konstruieren wir eine entartende "aquivariante quasi-projektive Familie, zu der  die gegebene Kurve geh\"ort und die gewisserma{\ss}en der einzige Baustein ist, aus dem sich jede gez\"ugelte \"aquivariante Familie, die \"uber eine feste Menge Punkte verzweigt, zusammensetzt.  
Das Ergebnis wird benutzt, um Automorphismen ordin\"arer Kurven induktiv zu untersuchen.}

\vspace{1cm}

\begin{center}
{\bf Degeneration of restrained group actions on curves} 
\end{center}

\noindent {\small {\bf Abstract.} 
An action of a finite group on a smooth projective curve over an algebraically closed field of positive characteristic is called {\sl restrained}, if all second ramification groups are trivial (e.g., every group action on an ordinary curve is restrained). When the ramification indices satisfy certain numerical criteria, we construct a degenerating equivariant quasi-projective family to which the given curve belongs, and which in a sense is the unique building block for all such restrained equivariant families that ramify above a fixed set of points. The result is used to inductively study automorphisms of ordinary curves.}

\vv

\vv

\newpage

{\bf Einf\"uhrung.} 

\vv

Es sei $X$ eine glatte projektive Kurve des Geschlechts $g$ \"uber einen algebraisch abgeschlossenen K\"orper der Charakteristik $p>0$, und $G$ eine endliche Gruppe Automorphismen von $X$. Sei $\pi : X \rightarrow Y= X / G$ die zugeh\"orige \"Uberlagerung. W\"ahrend im Falle komplexer Kurven alle Verzweigung zahm ist, d.h.\ die Isotropiegruppen von Punkten auf $X$ alle zyklisch sind, k\"onnen solche Isotropiegruppen \"uber K\"orper der positiven Charakteristik komplizierter aussehen: man nennt dieses Ph\"anomen {\sl wilde Verzweigung.} Es zeigt sich, dass die Verzweigungsgruppen $$G_x = G_{x,0} = \{ \sigma \in G : x^\sigma=x \}$$ eines Punktes $x \in X$ im allgemeinen Fall versehen sind mit einer absteigende Filtrierung durch Untergruppen, die sogenannten h\"oheren Verzweigungsgruppen
(cf.\ \cite{Ser}, IV):
$$G_{x,i}:=\{ \sigma \in G \ : \ \mbox{ord}_x(\sigma \pi_x-\pi_x)>i \} \ \ (i \geq 1),$$
wobei $\pi_x$ ein lokaler Uniformisierungsparameter ist in $x \in X$. 
Dann ist $G_{x,1}$ die einzige $p$-Sylowgruppe von $G_x$ und $G_x/G_{x,1}$ ist zyklisch von einer Ordnung teilerfremd zu $p$. 

Als n\"achste Stufe der Komplexit\"at nach dem zahmen Fall, d.h.\ $G_{x,1}=0$, schlagen wir eine Untersuchung der {\sl gez\"ugelten\footnote{\textrm{E.\ restrained, F.\ refr\'en\'e}} Verzweigung} vor, d.h.\ per Definition den Fall $G_{x,2}=0$. Dann ist also $G_x$ eine Gruppe
der Form $$G_x \cong ({\Z}/p)^{t_x} \sd \Z/n_x,$$ wobei $t_x$ und $n_x$ ganze Zahlen sind mit $n_x$ und $p$ teilerfremd. Wir nennen eine Gruppenoperation gez\"ugelt, falls die Verzweigung in allen Punkten $x \in X$ gez\"ugelt ist. Der Begriff ist aus verschiedenen Standpunkten nat\"urlich. So hat {\sc S.\ Nakajima} gezeigt (\cite{Nak}), dass jede Gruppenoperation auf einer ordin\"aren Kurve, d.h.\ eine Kurve, wof\"ur $\mbox{Jac}(X)[p]=(\Z/p)^g$, gez\"ugelt ist. Aber obwohl Ordinarit\"at eine generische Eigenschaft ist, d.h.\ Zariski-dicht im Modulraum, kann diese Eigenschaft trotzdem verloren gehen unter Entartung. Dagegen werden wir bald zeigen, dass die Eigenschaft $G_{x,2}=0$ erhalten bleiben {\sl kann}.

Der {\sc Grothendieck}sche Formalismus der \"aquivarianten Kohomologie (\cite{Gro}, \cite{BM}) erlaubt die vollst\"andige Berechnung der pro-darstellbaren Um\-h\"ullen\-den des infinitesimalen Deformationsfunktors von $(X,G)$ im Sinne {\sc Schlessinger}s (\cite{Sch}) als Funktion von $g$, der Anzahl der Verzweigungspunkte und der lokalen Verzweigungsdaten $(t_x,n_x)$ (cf.\ \cite{CK2}, Thm.\ 5.1): man bemerke, dass die Umh\"ullende nicht glatt zu sein braucht. 

Ziel dieser Arbeit ist die Erweiterung einer solchen infinitesimalen Richtung zu einer globalen quasi-projektiven Familie, deren Entartung voll\-st\"an\-dig verstanden werden kann. Es stellt sich heraus, dass im gez\"ugelten Falle die Operation von einem nicht-trivialen $G_{x,1}$ auf $k[[\pi_x]]$ durch M\"obius\-trans\-for\-mationen der Lokaluniformisierenden gegeben wird:
$$\pi_x \mapsto \frac{\pi_x}{u \pi_x +1}, u \in V, $$ wobei $V$ ein $t_x$-dimensionaler $\F_p$-Vektorraum innerhalb $k$ ist. Sei $s_x$ die multiplikative Ordnung von $p$ modulo $n_x$, d.h.\ $s_x=$ die Ordnung von $p$ in $(\Z/n)^\times$. Falls es $x \in X$ gibt mit $s_x \neq t_x$, dann kann man die Operation von $V$ so deformieren, dass ein ausgezeichneter Vektor $v$ in $V$ null wird: man l\"asst $v \mapsto 0$ in einer Zerlegung $V=V' \oplus \F_{p^{s_x}} \cdot v$. 

Diese Deformation l\"asst sich global fortsetzen (im Wesentlichen durch {\sc Grothendieck}s Algebraisierungssatz und das {\sc Artin}sche Approximationsverfahren), und es ergibt sich eine Familie von \"Uber\-la\-ge\-rungen \"uber eine feste Kurve, deren Verzweigungspunkte au{\ss}erdem fest sind --- der nat\"urliche Modulraum einer solchen \"Uber\-la\-ge\-rung ist affin, und also erlaubt die Familie eine entartete Faser (siehe {\sc Pries} \cite{Pri2} 3.3). Diese Entartung kann man
vollst\"andig verstehen, indem man sorgf\"altig die Moduln von $\F_p$-linearen Unterr\"aumen von $k$ und deren Fahnen untersucht, wobei die modulare Invariantentheorie von $\GL_n$ nat\"urlich ins Spiel kommt (\cite{Hew}). Die Hauptergebnisse werden im folgenden Satz beschrieben, der nicht nur besagt, dass eine Entartung existiert, sondern auch, dass die Operation auf jeder der horizontalen Komponenten einer solchen Entartung erneut gez\"ugelt ist, und au{\ss}erdem die Tr\"agheitsgruppen, d.h.\ Isotropiegruppen von Komponenten, beschreibt: 

\vv

{\bf Satz A.} {\sl Sei  $ \pi : X \rightarrow Y$ eine gez\"ugelte Galoissche \"Uberlagerung von Kurven \"uber einen algebraisch abgeschlossenen K\"orper $k$ der Charakteristik $p>0$ mit Galoisgruppe $G$, verzweigt \"uber $B \subset Y$, und es gebe einen Punkt $x \in X$ mit $1 \neq t_x \neq s_x$. Sei $y=\pi(x)$. Dann gibt es eine quasi-projektive glatte Kurve $S$, Punkte $\{0, \infty\} \in S$ und eine $G$-\"Uberlagerung von Kurven \"uber $S$:
\unitlength0.7pt
\begin{center}
\begin{picture}(100,45)
\put(15,35){$\Pi$}
\put(-10,25){$\X$}
\put(4,30){\vector(1,0){34}}
\put(40,25){$\mathcal{Y}$}
\put(0,20){\vector(1,-1){15}}
\put(42,20){\vector(-1,-1){15}}
\put(15,-10){$S$}
\end{picture}
\end{center}
mit

(i) $\Pi_\infty \cong \pi$;  

(ii) f\"ur alle $s \in S-\{0\}$ gilt $\mathcal{Y}_s = Y \times \{ s \}$ --- dabei ist $\Pi_s$ eine glatte $G$-\"Uberlagerung, verzweigt \"uber
$B \times \{ s \}$, derart dass die Verzweigung \"uber $(b,s)$ in $\Pi_s$ f\"ur $b \neq y$ isomorph ist zur Verzweigung \"uber $b$ in $\pi$;  

(iii) $\mathcal{Y}_{0}$ ist die Vereinigung von $Y$ und einer glatten rationalen      Kurve $E$, die $Y$ normal in $y$ schneidet; desweiteren ist $\Pi^{-1}(E)$ eine Vereinigung von glatten rationalen Kurven, die in je $p$ Punkten normal geschnitten werden von horizontalen Kurven, wobei die Einschr\"ankung der Quotientenabbildung auf solche horizontale Kurven gez\"ugelt ist;

(iv) jede zusammenh\"angende Komponente von $\Pi^{-1}(E) \subseteq \X_{0}$ hat Tr\"ag\-heits\-grad  $n_x p^{t_x-1}$, aber die horizontalen Kurven haben triviale Tr\"ag\-heits\-gruppe;

(v) Versieht man $\X_{S-\{0\}}$ mit einer gepunkteten Struktur, indem man alle Verzweigungspunkte auf jeder Faser $\X_s$ auszeichnet, so ist $\X_{S-\{0\}}$ eine stabile gepunktete Kurve, au{\ss}er im Fall wo $X=Y=\P^1, G=H$ und nur ein Punkt verzweigt ist; mit dieser Ausnahme \"ubertr\"agt sich die gepunktete Struktur nat\"urlich auf $\X$, wodurch $\X$ eine stabile gepunktete Familie von Kurven wird. } 

\vv

$$ \hspace{-2cm} \begin{picture}(0,0)%
\epsfig{file=special2.pstex}%
\end{picture}%
\setlength{\unitlength}{1973sp}%
\begingroup\makeatletter\ifx\SetFigFont\undefined%
\gdef\SetFigFont#1#2#3#4#5{%
  \reset@font\fontsize{#1}{#2pt}%
  \fontfamily{#3}\fontseries{#4}\fontshape{#5}%
  \selectfont}%
\fi\endgroup%
\begin{picture}(5925,8421)(601,-9199)
\put(3676,-8461){\makebox(0,0)[lb]{\smash{\SetFigFont{8}{9.6}{\familydefault}{\mddefault}{\updefault}$E$}}}
\put(6526,-7486){\makebox(0,0)[lb]{\smash{\SetFigFont{8}{9.6}{\familydefault}{\mddefault}{\updefault}$Y$}}}
\put(2401,-9136){\makebox(0,0)[lb]{\smash{\SetFigFont{8}{9.6}{\familydefault}{\mddefault}{\updefault}Die spezielle Faser $\X_0 \rightarrow \mathcal{Y}_0$}}}
\put(601,-7636){\makebox(0,0)[lb]{\smash{\SetFigFont{8}{9.6}{\familydefault}{\mddefault}{\updefault}$\mathcal{Y}_0$}}}
\put(3976,-7411){\makebox(0,0)[lb]{\smash{\SetFigFont{8}{9.6}{\familydefault}{\mddefault}{\updefault}y}}}
\put(4051,-6061){\makebox(0,0)[lb]{\smash{\SetFigFont{8}{9.6}{\familydefault}{\mddefault}{\updefault}$\Pi_0$}}}
\put(5251,-3361){\makebox(0,0)[lb]{\smash{\SetFigFont{8}{9.6}{\familydefault}{\mddefault}{\updefault}$\mathbf{P}^1$ mit Tr\"agheitsgrad $n_x p^{t_x-s_x}$}}}
\put(601,-2761){\makebox(0,0)[lb]{\smash{\SetFigFont{8}{9.6}{\familydefault}{\mddefault}{\updefault}$\X_0$}}}
\put(5926,-2161){\makebox(0,0)[lb]{\smash{\SetFigFont{8}{9.6}{\familydefault}{\mddefault}{\updefault} \} horizontale gez\"ugelte Kurven}}}
\put(5626,-4711){\makebox(0,0)[lb]{\smash{\SetFigFont{8}{9.6}{\familydefault}{\mddefault}{\updefault}  \} horizontale gez\"ugelte Kurven}}}
\end{picture}
 $$ 

\vv

\noindent Dieses Entartungsphenomen ist im folgenden Sinne sogar univer\-sell: w\"ah\-rend sich \"uber einem K\"orper der Charakteristik null bei Deformation einer \"Uber\-la\-ge\-rung zwangsl\"aufig die Verzweigungspunkte selbst bewegen, k\"onnen sie bei positiver Charakteristik fest bleiben. Aus der Berechnung des infinitesimalen versellen Deformationsrings in \cite{CK2} folgt, dass solche Deformationen nur unbehindert \"uber einen Ring der Krulldimension $>1$ stattfinden k\"onnen, falls $s_x \neq t_x$, und dann sind sie in der Tat durch Deformation der Einbettung $V \hookrightarrow k$ gegeben (\cite{CK2}, Abschnitt 4). Eine solche Einbettung entartet genau durch Verschwindenlassen eines Unterraumes $W$ von $V$. Wir haben im obigen Satz genau die grundlegende Situation beschrieben, wobei $W$ eindimensional ist, woraus sich der allgemeine Fall induktiv ableiten l\"asst.     

\vv

Das Ergebnis hat folgende Anwendung auf Automorphismengruppen gez\"ugelter Kurven, d.h.\ Kurven, f\"ur die $\pi_X: X \mapsto  X /\Aut(X)$ gez\"ugelt ist. Von {\sc S.\ Nakajima} (loc.\ cit.) stammt der Satz, dass eine gez\"ugelte Kurve $X$ des Geschlechts $g \geq 2$ h\"ochstens $|\Aut(X)| \leq 84 g (g-1)$ Automorphismen hat (der Beweis besteht im Grunde genommen aus geschickten Absch\"atzungen in der {\sc Riemann-Hurwitz-Zeuthen}schen Formel). Man kennt aber keine
Familie $\{X_i\}$, deren Geschlechter $g_i$ strikt steigend sind, f\"ur die $|\Aut(X_i)|$ ein Polynom vom Grad $3$ in $\sqrt{g_i}$ \"ubertrifft. In der Tat ist keine solche Familie bekannt mit
$|\Aut(X_i)| > \tilde{f}(g_i)$, wobei
$$ \tilde{f}(g) = \max \{ 84(g-1), 2 \sqrt{g} (\sqrt{g}-1)^2 \}. $$
Man wei{\ss}, dass f\"ur eine Mumfordkurve $X$ tats\"achlich $|\Aut(X)| \leq \tilde{f}(g)$ (\cite{CKK}); dieses Ergebniss ist scharf, und wir haben die Vermutung ge\"au{\ss}ert, es sei auch f\"ur ordin\"are Kurven richtig. 

Der obige Satz wird uns erlauben, die Zahlen $s_x$ induktiv zu erh\"ohen durch Deformation von $\pi_X$, bis die Gleichung $s_x=t_x$ erreicht wird. Nennen wir eine Kurve $X$ {\sl immobil}, wenn $ X / \Aut(X)  \cong \P^1$, $\pi_X$ gez\"ugelt ist, verzweigt \"uber zwei Punkten, oder \"uber drei Punkten mit Verzweigungstyp $(2,2,\cdot)$ und $n_x \neq 1$, aber $s_x=t_x$ f\"ur alle wilden Verzweigungspunkte $x \in X$. 

\vv

{\bf Satz B.} {\sl Sei $f : \Z_{\geq 2} \rightarrow \Z_{\geq 0}$ eine Funktion mit $f(g) \geq \tilde{f}(g)$ und $f(g)/(g-1)$ eine steigende Funktion von $g$. Unter der Annahme, dass f\"ur jede immobile Kurve $X$ von jedem Geschlecht $g \geq 2$ die Schranke $|\Aut(X)| \leq f(g)$ gelte, gilt eine solche Schranke auch f\"ur jede 
gez\"ugelte Kurve von jedem Geschlecht $g \geq 2$.}

\vv

Der Beweis ist im allgemeinen Fall induktiv und benutzt die besagte Existenz einer Entartung wie oben. Man bemerke aber, dass der Fall vollst\"andig multiplikativer Entartung einen anderen Beweis (mit Hilfe der {\sc Schottky}-Uniformisierung und einer Klassifikation Gruppenoperationen auf dem {\sc Bruhat-Tits}-Baum) erfordert, der im Grunde genommen in 
\cite{CKK} ausgef\"uhrt wurde, und ebenso der Fall, wo alle Komponenten elliptisch sind, wobei die Stabilit\"at der gepunkteten Kurve $\X_0$ ins Spiel kommt. 

\vv

Wir k\"onnen also unsere Vermutung zur\"uckf\"uhren auf die

\vv

{\bf Vermutung.} {\sl Jede immobile Kurve $X$ des Geschlechts $g \geq 2$ erf\"ullt $|\Aut(X)| \leq \tilde{f}(g)$.} 

\vv

\vv

Wir beschreiben schlie{\ss}lich kurz den Aufbau dieser Arbeit. In einem
ersten Abschnitt wird die lokale Theorie gez\"ugelter Gruppenoperationen untersucht, insbesondere wird anhand des {\sc Katz-Gabber}schen Lokal-Global-Prinzips eine kanonische Form f\"ur eine solche Operation \"uber komplette Bewertungsringe aufgestellt. Die Moduln solcher Operationen zeigen dadurch einen engen Zusammenhang mit der Variation
$\F_{p^{s_x}}$-linearer Unterr\"aume des K\"orpers $k$, der im n\"achsten Abschnitt
ausgearbeitet wird, anhand dessen dann schlie{\ss}lich  im Abschnitt \ref{familie} die f\"ur das Ergebnis erforderliche lokale Familie gez\"ugelter \"Uberlagerungen
konstruiert werden kann. Danach folgt in \ref{beweis} der Beweis des Hauptsatzes A und der Proposition. Anschlie{\ss}end wird im Abschnitt \ref{aut} Satz B gezeigt und kommentiert. Im ersten Anhang zeigen wir, dass die kanonische Form auch \"uber henselsche Ringe gilt. Im letzten Anhang
wird der Begriff "`Konfigurationsraum"' (aus {\sc Pries} \cite{Pri2}) erweitert zum Konzept "`Konfigurierender Funktor"', und wir zeigen anhand der Ergebnisse aus dem ersten Anhang, dass sich aus unserer Konstruktion im Abschnitt \ref{familie} ein nat\"urliches Beispiel f\"ur dieses verallgemeinerte
Konzept ergibt. 

\vv

\vv

\sectioning{Gez\"ugelte Verzweigung} \label{gez}

\vv

\paragraph\label{subpara-cov1} {\bf Galoistheoretische Vorbereitungen.} \ Wir stellen in den ersten Abschnitten einige Galoistheoretische Fakten bereit, die wir als solche nicht in der Literatur vorgefunden haben. Es sei $f:X \rightarrow Y$ ein Morphismus zwischen Schemata. Der Morphismus $f$ wird eine {\sl (endliche) \"Uberlagerung} (manchmal auch {\sl endliche Erweiterung}) genannt, falls $f$ endlich und unverzweigt (im Sinne von \cite{EGA}, $\textrm{IV}_4$ (17.3.7)) ist in jedem Punkt der H\"ohe Null auf $X$. Wenn $f : X \rightarrow Y$ eine solche \"Uberlagerung ist und $Y$ lokal noethersch und reduziert, dann gibt es eine dichte offene Teilmenge $U \subseteq Y$, \"uber der $f$ flach ist (\cite{EGA}, $\textrm{IV}_2$ (6.9.1)); in Folge dessen ist $f$ generisch \'etale. 
Wir schreiben $G = \Aut_Y(X)$ und nehmen ab jetzt an, die Gruppe operiert von rechts auf $X$ und ist endlich. Letztere Eigenschaft ist automatisch erf\"ullt, falls $X$ endlich viele irreduzibele Komponenten hat. Weil $f$ affin ist, macht es Sinn, den Quotienten $X/G$ zu bilden (\cite{SGA1} Expos\'e V, Prop.\ 1.1). Falls die von $f$ induzierte Abbildung $X/G \rightarrow Y$ ein
Isomorphismus ist, nennen wir $f$ eine {\sl Galoissche \"Uberlagerung} (oder {\sl Galoiserweiterung}) {\sl mit Galoisgruppe $G$}, und schreiben wir $\Gal(X/Y)=G$. 

\vv

\paragraph\label{subpara-cov4} {\bf Isomorphismus von Galoiserweiterungen.} \ Ein Isomorphismus zwischen zwei  Galoisschen \"Uberlagerungen $f: X \rightarrow Y$ und $f' : X' \rightarrow Y$ ist ein Paar $(\phi,h)$, wobei $\phi: X \rightarrow X'$ ein Isomorphismus ist zwischen $X$ und $X'$ \"uber $Y$ und $h: \Gal(X/Y) \stackrel{\sim}{\rightarrow} \Gal(X'/Y)$ ein Gruppenisomorphismus, derart dass $h(\sigma) = \phi \circ \sigma \circ \phi^{-1}$ f\"ur alle $\sigma \in \Gal(X/Y)$. Die Beweise der n\"achsten zwei Lemmata lassen wir ihrer Einfachkeit halber aus.

\vv

\paragraph\label{lem-covaff1}
 {\bf Lemma.} \ {\sl Sei $B$ ein ganzabgeschlossener Ring mit Quotientenk\"orper $L$, und sei $G$ eine endliche Gruppe, die von links auf $B$ operiert via Ringautomorphismen. Bezeichne mit $A$ den Invariantenring $A=B^G$. Dann gilt:

(i) Der Fixk\"orper $L^G$ ist Quotientenk\"orper von $A$;

(ii) $B/A$ ist eine Galoiserweiterung mit Galoisgruppe $G$;

(iii) $A$ ist ganzabgeschlossen; au{\ss}erdem ist $A$ ein (diskreter) Bewertungsring, falls $B$ es ist. \qed}

\vv

\paragraph\label{lem-covaff2}
 {\bf Lemma.} \ {\sl Sei $L$ eine endliche Galoiserweiterung des Quotienten\-k\"or\-pers $K$ eines ganzabgeschlossenen Ringes $A$, und $B$ der ganze Abschlu{\ss} von $A$ in $L$. Dann ist $B/A$ eine Galoiserweiterung mit Galoisgruppe $\Gal(L/K)$. \qed }

\vv

\paragraph\label{subpara-norcov1} {\bf Die Galoiskorrespondenz.} \ Sei $f : X\rightarrow Y$ eine Galoissche \"Uberlagerung von lokal noetherschen normalen ganzen Sche\-ma\-ta und $Z$ ein zwischen $X$ und $Y$ liegendes normales ganzes Schema, derart dass $Z \rightarrow Y$ eine \"Uberlagerung ist. Dann ist $X \rightarrow Z$ auch eine \"Uberlagerung (\cite{EGA}, II (6.1.5) (v) und $\textrm{IV}_4$ (17.3.4)), und also ist $X$ die Normalisierung von $Z$ im Funktionenk\"orper $\kappa(X)$ von $X$. Aus (\ref{lem-covaff2}) folgt, dass $X/Z$ eine Galoissche \"Uberlagerung mit Galoisgruppe $H_Z:=\Gal(\kappa(X)/\kappa(Z)) = \Aut_Z(X)$ ist. Eine alternative Definition von $H_Z$ lautet: f\"ur jedes $\sigma \in G=\Gal(X/Y)$ hat man die konjugierte \"Uberlagerung $Z^\sigma \rightarrow Y$, wobei es einen 
Isomorphismus $\sigma_Z : Z \stackrel{\sim}{\rightarrow} Z^\sigma$ gibt. Dann ist $H_Z$ die Untergruppe der Elemente $\sigma$ von $G$, f\"ur die das dazugeh\"orige $\sigma_Z$ mit den \"Uberlagerungsabbildungen $X \rightarrow Z$ und $X \rightarrow Z^\sigma$ kommutiert. Die Korrespondenz $Z \leadsto H_Z$ induziert also eine Inklusionsumkehrende Bijektion zwischen Isomorphismusklassen von Zwischen\"uberlagerungen von $f$ und Untergruppen von $G$, wobei Galoissche
Zwischen\"uberlagerungen \"uber $Y$ mit normalen Untergruppen korrespondieren. 

\vv

\paragraph\label{nolabel} {\bf Bemerkung.} \  In der Kategorie der quasiprojektiven Kurven \"uber einen algebraisch abgeschlossenen K\"orper haben wir soeben die Galoiskorrespondenz {\sl glatter} Kurven geschildert, wegen (\cite{Sai}, Prop.\ 1.6) l\"a{\ss}t diese sich aber auch auf allen {\sl semistabilen} Kurven mit glatten irreduzibelen Komponenten fortsetzen.

\vv

\paragraph\label{subpara-ramgroup}
 {\bf Verzweigungsgruppen.} \ Sei $k$ ein K\"orper der Charakteristik $p>0$ und $A$ ein diskreter Bewertungsring \"uber $k$ mit Restklassenk\"orper $k$; sei $B/A$ eine endliche Galoiserweiterung mit Galoisgruppe $G=\Gal(B/A)$. 
Wir bezeichnen mit  $K$ und $L$ die Quotientenk\"orper von $A$ und $B$, und mit $\m_A,\m_B$ deren maximale Ideale. 

Der Ring $B$ ist auch diskret bewertet und $\Spec B \rightarrow \Spec A$ ist endlich und flach (\cite{Mat}, 23.1), insbesondere ist $B$ ein freies $A$-Modul endlichen Ranges. Sei $x$  eine Uniformisierende aus $B$ und $v_B : B \rightarrow \Z \cup \{ \infty \}$ die dazugeh\"orige normalisierte Bewertung. F\"ur $i \geq -1$ definieren wir die $i$-te Verzweigungsgruppe $G_i$ durch
$$ G_i := \{ \sigma \in G : v_B(\sigma(x) - x) > i \}. $$
Die $G_i$ bilden eine absteigende Filtrierung der Gruppe $G$ durch normale Untergruppen: 
$$
G=G_{-1}\supseteq G_0\supseteq G_1\supseteq\cdots .
$$  

\vv

\paragraph\label{dfn-restrained} {\bf Gez\"ugelte Erweiterungen.} \ Wir nennen eine Galoiserweiterung $B/A$ {\sl gez\"ugelt}, falls $B/A$ total verzweigt ist und $G_2=0$. 

\vv

\paragraph\label{rem-restrained} {\bf Bemerkung.} \ Die Totalverzweigung von $B/A$ bedeutet, 
dass die induzierte Rest\-klassen\-k\"or\-per\-er\-wei\-te\-rung trivial ist; insbesondere ist $G=G_0$. Im Falle eines algebraisch abgeschlossenen Restklassenk\"orpers ist diese Bedingung also immer erf\"ullt.  

\vv

\paragraph\label{prop-subquot}
 {\bf Proposition.} \ {\sl Wenn $B/A$ eine gez\"ugelte Galoiserweiterung mit Galoisgruppe $G$ ist und $H$ eine Untergruppe von $G$, so ist die Galoiserweiterung $B/B^H$ erneut gez\"ugelt. Ist au{\ss}erdem $H$ eine normale Untergruppe, so ist auch $B^H/A$ gez\"ugelt.}

\vv

\pf Die erste Aussage ist klar. Zum Beweis der zweiten Aussage nehme man an, $B^H/A$ sei nicht gez\"ugelt. Da die Erweiterung offensichtlich total verzweigt ist, bedeutet dies $(G/H)_2 \neq 0$. Die Herbrandsche Funktion $$\phi_{L/L^H}(u):=\int_0^u \frac{dt}{[G_0:G_t]}$$ erf\"ullt also 
$ \phi_{L/L^H}(u)>1 \Rightarrow G_u \not \subset H $
 (\cite{Ser}, IV, Lemme 5). Deswegen muss $u \leq 1$, und also 
$$ 1 <\varphi_{L/L^H}(u)=u\frac{|H_1|}{|H_0|}\leq 1,$$
ein Wiederspruch. \qed
 
\vv

\paragraph\label{Nak-thm}
 {\bf Satz} ({\sc S.\ Nakajima}, \cite{Nak} 2(i)). \ {\sl Sei $X/k$ eine glatte projektive Kurve vom Geschlecht $g>0$ \"uber einen algebraisch abgeschlossenen K\"orper $k$ der Charakteristik $p>0$, und nehme an, $X$ ist ordin\"ar, d.h. der $p$-Rang der Jacobischen $\textrm{Jac}(X)$ ist maximal: $\textrm{Jac}(X)[p] = (\Z/p)^g$. Dann ist f\"ur jede endliche Gruppe von Automorphismen $G \subseteq \Aut_k(X)$ und jeden geschlossenen Punkt $x \in X$, die Erweiterung  $\O_{X,x}/\O^{G_x}_{X,x}$ von lokalen Ringen gez\"ugelt, wobei $G_x$ die Zerlegungsgruppe von $x$ ist.} 

\vv

\paragraph\label{Nak-rm} 
{\bf Bemerkung.} \ Generische Kurven sind ordin\"ar, d.h.\ Ordinarit\"at ist Zariski-dicht im Modulraum von Kurven vom Geschlecht $g>0$ (vgl.\ \cite{Oort}). 

\vv

\paragraph\label{subpara-group1} {\bf Die Gruppe $G_{n,d}$.} \ Seien $n$ und $d$ zwei ganze Zahlen mit $n \geq 1$, $n$ teilerfremd zu $p$ und $d \geq 0$. Wir setzen $q=p^s$, wobei $s$ die kleinste positive ganze Zahl ist, f\"ur die $n$ die Zahl $p^s-1$ teilt --- was gleichwertig damit ist, dass $s=[\F_p(\zeta):\F_p]$ mit $\zeta$ einer festgew\"ahlten primitiven $n$-ten Einheitswurzel in $k$; oder auch damit, dass $p$ die Ordnung $s$ hat in der multiplikativen Gruppe $(\Z/n)^\times$ der ganzen Zahlen modulo $n$. Wir bezeichnen mit $G_{n,d}$ die abstrakte Gruppe $G_{n,d} = H \sd \Z/n$, wobei $H=(\F_q^d,+)$ und das semidirekte Produkt durch den Morphismus $\phi : \Z/n \rightarrow \Aut H : a \mapsto (\sigma \mapsto \zeta^a \sigma)$ gegeben wird, d.h.\ in expliziter Form
$$
(\sigma,a)(\tau,b)=(\sigma+\zeta^a\tau,a+b),
$$
mit $\sigma,\tau\in H$ und $a,b\in\Z/n\Z$.
Wir nennen {\sl Torusteil} von $G_{n,d}$ jede Teilgruppe der Ordnung $n$. 

\vv

\paragraph\label{lem-group2}
{\bf Lemma.}\ {\sl (i) Die Gruppe $H$ ist die einzige $p$-Sylowgruppe von $G_{n,d}$, und $H$ ist elementar abelsch der Ordnung $q^d$.

(ii) Jeder Torusteil von $G_{n,d}$ ist zyklisch und hat ein erzeugendes Element  $\gamma$ mit $\gamma \sigma \gamma^{-1} = \zeta \sigma$ f\"ur $\sigma \in H$. 

(iii) Seien $T$ und $T'$ zwei Torusteile von $G_{n,d}$. Dann gibt es $\sigma \in H$ mit $T'=\sigma T \sigma^{-1}$.}

\vv

\pf Da {\sl (i)} und {\sl (ii)} offensichtlich sind, brauchen wir nur {\sl (iii)} zu zeigen. Sei $\gamma$ ein erzeugendes Element f\"ur $T$ mit $\gamma \sigma \gamma^{-1} = \zeta \sigma$ f\"ur $\sigma \in H$. Dann gibt es $\tau \in H$, so dass $\gamma' = \tau \gamma$ die Gruppe $T'$ erzeugt. Man w\"ahle $\sigma = (1-\zeta)^{-1} \tau$. \qed

\vv

\paragraph\label{dfn-type} {\bf Gez\"ugelte Erweiterungen vom Typ $(n,d)$.} \ Sei $B/A$ eine gez\"ugelte Erweiterung mit Galoisgruppe $G$ wie oben. Aus (\cite{Ser} IV, Cor.\ 1, Cor.\ 3, Cor.\ 4, Prop.\ 9) ergibt sich direkt, dass $G$ zu einer Gruppe $G_{n,d}$ isomorph ist. Wir nennen $B/A$ dann eine gez\"ugelte Erweiterung {\sl vom Typ $(n,d)$}. 

\vv

\paragraph\label{subpara-map1} {\bf Die Abbildung $\theta_{B/A}$.} \ Sei $B/A$ vom Typ $(n,d)$ und $H=G_1$ die $p$-Sylowgruppe von $G$. Sei $\Omega_{B/A} = \m_B/\m_B^2$ der Kotangentialraum von $B$ und $x$ eine Uniformisierende von $B$. Wie in (\cite{Ser}, IV.2) definiert man eine $\F_p$-lineare Einbettung 
$$
\theta_{B/A}\colon H\longrightarrow\Omega_B,\qquad
\sigma\mapsto\frac{\sigma(x)-x}{x}\ \textrm{mod}\ \m^2_B
$$
(dort $\theta_1$ genannt), die von der Wahl von $x$ unabh\"angig ist. Aus (\cite{Ser} IV, Prop.\ 9) folgt unmittelbar, dass $\theta_{B/A}$ mit der Multiplikation mit $\zeta$ kommutiert. Deshalb ist $\theta_{B/A}$ nicht nur $\F_p$-, sondern sogar $\F_q$-linear. 

\vv

\paragraph\label{subpara-map2} {\bf Die Invariante $\omega_{B/A}$.} \ Es ist $\theta:=\theta_{B/A}/dx$ eine $\F_q$-lineare Einbettung von $H$ in $k$. Entwickeln wir f\"ur $\sigma \in H$ das Element $\sigma(x)$ als Potenzreihe in $x$ innerhalb der Vervollst\"andigung $\hat{B}=k[[x]]$ von $B$, dann finden wir
$$
\sigma(x)=x+\theta(\sigma)x^2+\cdots.
$$
Man bemerke, das $\theta$ von der Wahl von $x$ abh\"angt, aber nur durch Multiplikation mit einem Skalar aus $k$. Wir definieren die {\sl $\omega$-Invariante} $\omega_{B/A}$ von $B/A$ als die Klasse modulo $k^\times$ des $\F_q$-linearen Unterraums $\theta(H) \subset k$. 

\vv

\paragraph \label{prop-canonicalform} {\bf Proposition} (Kanonische Form gez\"ugelter Gruppenoperationen). \ {\sl Sei $A$ ein vollst\"andiger diskreter Bewertungsring und $B/A$ eine gez\"ugelte Erweiterung vom Typ $(n,d)$ mit Galoisgruppe $G$ mit $p$-Sylowgruppe $H$. W\"ahle einen festen Torusteil $T \subset G$. Dann gibt es ein Uniformisierende $x$ f\"ur $B$ und ein
erzeugendes Element $\gamma$ von $T$, derart dass die Operation von $G$ gegeben ist durch:
$$
\gamma(x)=\zeta x,\qquad\sigma(x)=\frac{x}{1-\theta(\sigma)x} \mbox{ \ \ \ f\"ur jedes $\sigma\in H$,}
$$
wobei $\theta=\theta_{B/A}/dy$.}

\vv

\pf Aus dem {\sc Katz-Gabber}schen Lokal-Global-Prinzip (\cite{Katz}, 1.4)  ergibt sich, dass die Galoissche \"Uberlagerung $\Spec B \rightarrow \Spec A$ die Spezialisierung einer speziellen $G$-\"Uberlagerung von Kurven $X \rightarrow \P^1_k$ \"uber $k$ im $\infty$ ist, die nur \"uber $0$ und $\infty$ verzweigt ist. Weil die Verzweigung \"uber 0 zwangsl\"aufig zahm ist und keine $\Z/n$-\"Uberlagerung der projektiven Gerade in nur einem Punkt verzweigt sein kann, ist die Zerlegungsgruppe \"uber $0$ gleich $\Z/n$. Da die Zerlegungsgruppe von $\infty$ per Konstruktion $G$ ist, ergibt sich aus der {\sc Riemann-Hurwitz-Zeuthen}schen Formel (\cite{Sti} 1.1) unmittelbar $X \cong \P^1_k$, weswegen $G$ durch M\"obiustransformationen operiert. 

Sei $x$ eine inhomogene Koordinate auf $X$, so dass $x=0$ \"uber $\infty$ liegt, und $x=\infty$ ein Punkt \"uber $0$ ist, der von $T$ festgehalten wird. Da $T$ die Fixpunkte $0$ und $\infty$ hat, gibt es ein erzeugendes Element $\gamma$ von $T$ mit $\gamma(x)=\zeta x$. 

Da $\sigma \in H$ nur $\infty$ fixiert, ist $$\sigma(x) = \frac{ax}{b+cx}$$ f\"ur $a,c \in k^\times$ und $b \in k$. Weil $\sigma(x)=x \mbox{ mod } \m^2_B$ in $B=\hat{\O}_{X,0}$, ist $a=b$ und wir k\"onnen $$\sigma(x) = \frac{x}{1+cx}$$ voraussetzen. Dann ist $$\theta_{B/A}(\sigma) = \frac{\sigma(x)-x}{x} \mbox{ mod } \m_B^2 = -c dy,$$ weswegen $c=-\theta(\sigma)$. \qed

\vv

\paragraph \ {\bf Bemerkungen.} \ (i) Wenn $H \neq 0$, d.h.\ $d>0$, so ist das  erzeugende Element $\gamma$ von $T$ eindeutig durch die Beziehung $\gamma \sigma \gamma^{-1} = \zeta \sigma$ f\"ur $\sigma \in H$ bestimmt. 

(ii) Statt $A$ vollst\"andig gen\"ugt es, $A$ henselsch vorauszusetzen: dies wird sich sp\"ater daraus ergeben, dass die Uniformisierende $x$ in der Proposition genau eine L\"osung des Polynoms $Q$ wie nachstehend in (\ref{prop-invariant}) ist, weswegen es approximiert werden kann. Wir verschieben den Beweis auf Anhang \ref{hens-norm}. 

(iii) Die Reduktion auf eine kanonische Form  aus \ref{prop-canonicalform} wurde in \cite{CK2}, 3.1 stillschweigend vorausgesetzt.

\vv

\paragraph\label{prop-isomaction}
{\bf Proposition.}\ {\sl Es sei $A$ vollst\"andig. Dann sind zwei gez\"ugelte Galoiserweiterungen $B/A$ und $B'/A$  vom Typ $(n,d)$ isomorph im Sinne von (\ref{subpara-cov4}) genau dann, wenn ihre $\omega$-Invarianten \"ubereinstimmen: $\omega_{B/A} = \omega_{B'/A}$.}

\vv

\pf Es operiere $G=G_{n,d}$ auf $B$ und $B'$ \"uber $A$, und es sei $H$ die $p$-Sylowgruppe und $T$ ein festgew\"ahlte Torusteil von $G$ mit erzeugendem Element $\gamma$. Wir schreiben die Operation von $\alpha \in G$ auf $B$ wie $\alpha$, aber auf $B'$ wie $\alpha'$. 

Wir k\"onnen  wie im Beweis von (\ref{prop-canonicalform}) voraussetzen, dass 
$B=\hat{\mathcal{O}}_{P,\xi}, B'=\hat{\mathcal{O}}_{P',\xi'}$ und $A=\hat{\mathcal{O}}_{P'',\eta}$ mit zugeh\"origen Gruppenoperationen, wobei
$P, P'$ und $P''$ drei projektive Geraden sind mit den jeweiligen inhomogenen Koordinaten $x$ und $x'$; $\xi$ (resp.\ $\xi'$) der Punkt $\{x=0\}$ (resp.\ $\{ x'=0 \}$) ist; G auf $P$ und $P'$ operiert durch  
$$ \left\{ \begin{array}{l}
\gamma(x)=\zeta x,\\ {\displaystyle \sigma(x)=\frac{x}{1-\theta(\sigma)x} }, 
\end{array}\right. \ \ \left\{ \begin{array}{l}
\gamma'(x')=\zeta x',\\ {\displaystyle \sigma'(x')=\frac{x'}{1-\theta'(\sigma)x'} } \end{array}\right. $$
f\"ur jedes $\sigma\in H$, wobei $\theta$ und $\theta'$ Einbettungen von $H$ in $k$ sind; $\gamma \sigma \gamma^{-1} = \zeta \sigma$ gilt, und $\pi \colon P \rightarrow P''$ und $\pi' \colon P' \rightarrow P''$ die zugeh\"origen  Quotientenabbildungen sind, die beide \"uber den gleichen Punkt $\eta \in P''$ wildverzweigt sind. Insbesondere ist $\theta=\theta_{B/A}/dx$ und $\theta'=\theta_{B'/A}/dx'$.  

Nehmen wir zuerst an, $B/A$ und $B'/A$ seien isomorphe $G$-Erweiterungen, und sei $(\phi,h)$ das zugeh\"orige Paar mit $\phi: B \stackrel{\sim}{\rightarrow} B'$ der Isomorphismus \"uber $A$ und $h \in \Aut(G)$ ein Automorphismus 
mit $\alpha'= \phi \circ h(\alpha) \circ \phi^{-1}$. Weil der Urbildfunktor aus (\cite{Katz} 1.4.1) voll treu ist, sind die zugeh\"origen Operationen von $G$  auf $P$ und $P'$ isomorph. Deswegen k\"onnen wir $\phi$ als M\"obiustransformation $$\phi(x)=\frac{ax'+b}{cx'+d}$$ schreiben, mit $a,b,c,d \in k$. Weil $\sigma' \circ \phi = \phi \circ h(\sigma)$ f\"ur $\sigma \in H$ damit bedeutet, dass 
$$
\frac{(a-\theta'(\sigma)b)x'+b}{(c-\theta'(\sigma)d)x'+d}=
\frac{ax'+b}{(c-\theta(h(\sigma))a)x'+(d-\theta(h(\sigma))b)},
$$
folgt $b=0$ und $\theta'(\sigma)d=\theta(h(\sigma))a$. Letztere Gleichung impliziert, dass die Bilder von $\theta$ und $\theta'$ sich nur um einen Skalar in $k^\times$ unterscheiden, weswegen $\omega_{B/A}=\omega_{B'/A}$. 

Gilt umgekehrt $\omega_{B/A}=\omega_{B'/A}$, dann w\"ahlen wir $c \in k^\times$ mit $c\theta(H) = \theta'(H)$. Es ergibt sich 
$$
\sigma'(cx')=\frac{cx'}{1-(c^{-1}\theta'(\sigma))cx'}.
$$
Wenn wir also $x'$ durch $cx'$ ersetzen, k\"onnen wir annehmen, dass $\theta(H) = \theta(H')$. Es gibt also einen $\F_q$-linearen Automorphismus $h \in \Aut(H)$ mit $\theta'(\sigma)=\theta(h(\sigma))$, der sich, weil $h$ linear ist und die Gruppe explizit wie \ref{subpara-group1} beschrieben wird, durch $h(\gamma)=\gamma$ erweitern l\"asst zu einem Automorphismus von $G$. Dann ist
$(\phi,h)$ mit $\phi \ \colon \ B=k[[x]] \rightarrow B'=k[[x']]\ \colon \ x \mapsto x$ der erw\"unschte Isomorphismus der Operationen. \qed

\vv

Die Struktur einer gez\"ugelten \"Uberlagerung wird also vollst\"andig von der $\omega$-Invariante bestimmt. Da sie einen $\F_q$-linearen Unterraum vom K\"orper $k$ bildet, werden wir uns im n\"achsten Abschnitt zuerst dem Studium von Moduln solcher R\"aume widmen. 

\vv

\vv

\sectioning{Moduln $\F_p$-linearer Unterr\"aume eines K\"orpers}\label{unter}

\vv

Es sei $k$ ein K\"orper der Charakteristik $p>0$ und $q=p^n$. Die Modultheorie linearer Unterr\"aume von $k$ ist eng mit der klassischen {\sc Dickson}schen 
Invariantentheorie verkn\"upft. 

\vv

\paragraph\label{subpara-linearpoly1}
{\bf $\F_q$-lineare Polynome.} \ (i) Sei $R$ eine $\F_q$-Algebra. Wir nennen einen $\F_q$-linearen Unterraum $V \subseteq R$ {\sl nicht-entartet}, falls $V \cap \mathfrak{p} = \{ 0 \}$ f\"ur jedes Primideal $\mathfrak{p}$ von $R$, d.h.\ dass die Dimension von $V$ sich durch keine Spezialisierung \"andert; dies ist gleichbedeutend damit, dass $V-\{0\} \subseteq R^\times$. 

(ii) Ein Polynom $P$ in $R[X]$ der Form 
$$
P(a_0,a_1,\ldots,a_d)(X)=a_dX^{q^d}+a_{d-1}X^{q^{d-1}}+\cdots+a_1X^q+a_0X
$$
nennen wir ein {\sl $\F_q$-lineares Polynom der H\"ohe $d$ \"uber $\F_q$}. 

(iii) Ein Polynom wie in (ii) hei{\ss}t {\sl nicht-entartet}, wenn $a_0,a_d \in R^\times$. 

(iv) Ein nicht-entartetes $\F_q$-lineares Polynom $P$ {\sl zerf\"allt in $R$}, falls es einen nicht-entarteten $\F_q$-linearen Unterraum $V$ von $R$ und $c \in R^\times$ gibt mit $$P= c \cdot \prod_{v \in V}(X-v).$$ 

\vv

\paragraph\label{lem-ring1}
{\bf Lemma.}\ {\sl Sei $R$ eine $\F_q$-Algebra. Dann ist die Abbildung
\begin{center}
\begin{tabular}{rcl}
$\left\{
\begin{minipage}{8em}
\setlength{\baselineskip}{.85\baselineskip}
\begin{center}\begin{small}
{\slshape
nicht-entarteter $\F_q$-linearer Unterraum $V$ der Dimension $d$ in $R$}
\end{small}
\end{center}
\end{minipage}
\right\}$ &
$\stackrel{F}{\longrightarrow}$ &
$\left\{
\begin{minipage}{8em}
\setlength{\baselineskip}{.85\baselineskip}
\begin{center}\begin{small}
{\slshape
nicht-entartetes $\F_q$-lineares Polynom vom Grad $q^d$, das zerf\"allt in $R$}
\end{small}
\end{center}
\end{minipage}
\right\}\big/R^{\times}$ \\
& & \\
$V$ & $\mapsto$ & $\displaystyle \prod_{v\in V}(X-v)$
\end{tabular}
\end{center}
bijektiv mit Inverse
$P \mapsto \{ v \in R^\times \ : \ P(v)=0 \} \cup \{ 0 \}$.
}

\vv

\pf Sei $R$ zun\"achst ganz. Man kann $R$ in einen algebraisch abgeschlossenen K\"orper einbetten, und aus \cite{Gos}, 1.2.1 und 1.2.2 ergibt sich sofort die Bijektion
\begin{center}
\begin{tabular}{rcl}
$\left\{
\begin{minipage}{10.5em}
\setlength{\baselineskip}{.85\baselineskip}
\begin{center}\begin{small}
{\slshape
$\F_q$-linearer Unterraum $V$ von $R$ der Dimension $d$}
\end{small}
\end{center}
\end{minipage}
\right\}$ &
$\longrightarrow$ &
$\left\{
\begin{minipage}{10em}
\setlength{\baselineskip}{.85\baselineskip}
\begin{center}
\begin{small}
{\slshape
$\F_q$-lineares Polynom von Grad $q^d$ alle dessen Wurzeln in $R$ liegen
}
\end{small}
\end{center}
\end{minipage}
\right\}\big/R^{\times}$ \\

& & \\
$V$ & $\mapsto$ & $\displaystyle \prod_{v\in V}(X-v).$ \end{tabular} 
\end{center}
Diese Bijektion bildet offenbar nicht-entartete $\F_q$-lineare Unterr\"aume genau auf nicht-entartete $\F_q$-lineare Polynome ab: sie ist kompatibel mit jeder Spezialisierung $R \rightarrow R/\mathfrak{p}$ f\"ur ein Primideal $\mathfrak{p}$. 

Es sei jetzt $R$ eine $\F_q$-Algebra, nicht notwendig ganz. Dann ist $F$ offenbar immer noch surjektiv. Die Injektivit\"at ergibt sich aus dem Sachverhalt:
\begin{quote}
Sei $P=\prod_{v \in V} (X-v)$ ein nicht-entartetes $\F_q$-lineares Polynom vom Grad $q^d$, das in $R$ zerf\"allt, und $V$ sei ein nicht-entarteter $\F_q$-linearer Unterraum von $R$ der Dimension $d$. Sei $W=\{ w \in R^\times : \ P(w)=0 \} \cup \{ 0 \}$. Dann ist in der Tat $W=V$. 
\end{quote}
Offenbar ist $V \subseteq W$. Sei $\mathfrak p$ ein Primideal von $R$ und betrachte die kanonische Abbildung $\pi \ : \ R \rightarrow R/\mathfrak{p}$. Da $V$ nicht-entartet ist, gilt $\ker(\pi|_V) = V \cap \mathfrak{p} = 0$, also $V \cong \pi(V)$. Da $W \subseteq R^\times \cup \{ 0\}$, ist $W \cap \mathfrak{p} = 0$, weswegen $W \cong \pi(W)$. Weil das induzierte Polynom $\pi(P)$ genau $\pi(V)$ als Menge seiner Wurzeln hat, ist $\pi(W) \cong \pi(V)$. Wir schliessen hieraus, dass $\dim W = \dim V$, und da $V \subseteq W$, also $V=W$. \qed

\vv

\paragraph\label{subpara-GLquot1def} {\bf Quotient durch $\GL(d,\F_q)$.} Man l\"asst die Gruppe $\mathfrak{G}=\GL(d,\F_q)$ von rechts auf dem Polynomring $\F_q[U_1,\ldots,U_d]$ operieren durch $U_j\cdot A=\sum_{i=1}^da_{ji}U_i$ wobei $A=(a_{ij})$. Sei $\til{V} = \sum_{i=1}^d \F_q \cdot U_i$, und  $$T_0,\ldots,T_{d-1}\in\F_q[U_1,\ldots,U_d]$$ definiert durch die Gleichung 
$$
\prod_{v\in\til{V}}(X-v)=P(T_0,\ldots,T_{d-1},1)(X).
$$
Dann sind $T_i$ deutlich $\mathfrak{G}$-Invarianten.

\vv

\paragraph \label{lem-GLquot1}
{\bf Proposition} ({\sc Dickson} \cite{Dic}).\ {\sl Die Elemente $T_0,\ldots,T_{d-1}$ sind algebraisch unabh\"angig \"uber $\F_q$, und erzeugen den Invariantenring $$\F_q[U_1,\ldots,U_d]^{\mathfrak{G}}=\F_q[T_0,\ldots,T_{d-1}].$$
Insbesondere ist $\F_q[T_0,\ldots,T_{d-1}]$ ganzabgeschlossen im Zerf\"allungsk\"orper von $$P(T_0,\ldots,T_{d-1},1).$$}
\pf Siehe z.B.\ \cite{Wil}. \qed

\vv

\paragraph\label{cor-GLquot1}
{\bf Korollar.}\ {\sl Bezeichnet man mit $W$ das Produkt aller nicht-Null Elemente aus $\til{V}=\sum^d_{i=1}
\F_qU_i$, so ist 
$$
\F_q[U_1,\ldots,U_d,W^{-1}]^{\mathfrak{G}}=\F_q[T'_1,\ldots,T'_d,T^{\prime -1}_d],
$$
mit $P(1,T'_1,\ldots,T'_d)=Y\prod_{v\in\til{V}-\{0\}}(1-v^{-1}Y)$.}
\begin{eqnarray*}
\mbox{ {\sl Beweis.} \ } \F_q[U_1,\ldots,U_d,W^{-1}]^{\mathfrak{G}} &=& \F_q[T_0,\ldots,T_{d-1},T^{-1}_0] \\ &=&\F_q[T_1T^{-1}_0,\ldots,
T_{d-1}T^{-1}_0,T^{-1}_0,T_0] \\ &=& \F_q[T'_1,\ldots,T'_d,T^{\prime -1}_d]. \ \ \ \ \ \ \ \ \  \square
\end{eqnarray*}

\paragraph \label{subpara-functor1} {\bf Moduln von Unterr\"aumen.} \
Sei $k$ ein K\"orper der $\F_q$ umfasst, und betrachte den kovarianten Funktor 
$\mathbf{L}_d\colon\mathbf{Alg}_k\rightarrow\mathbf{Mengen},$
der einer $k$-Algebra $R$ die Menge $\mathbf{L}_d(R)$ aller nicht-entartete $\F_q$-linearen
Unterr\"aume der Dimension $d$ in $R$ zuordnet.
Die \'etale Garbifizierung dieses Funktors ist der kontravariante Funktor
$\mathcal{L}_d\colon\mathbf{Sch}_k\rightarrow\mathbf{Mengen}$  
$$
X\mapsto\mathcal{L}_d(X)=
\left\{
\begin{minipage}{17em}
\setlength{\baselineskip}{.85\baselineskip}
\begin{small} \begin{center}
{\slshape $\F_q$-Untermoduln von $\O^{\mbox{\tiny \'et}}_X$, die ein \'etales \\ $\F_q$-lokales System vom Rang $d$ bilden}
\end{center} \end{small}
\end{minipage}
\right\}.
$$
Wir ben\"otigen sp\"ater zus\"atzlich noch folgende Konstruktion. Sei $H=\F^d_q$ versehen mit der Standardbasis, und definiere den Funktor $\til{\mathbf{L}}_d\colon\mathbf{Alg}_k\rightarrow\mathbf{Mengen}$, der jeder 
$k$-Algebra $R$ die Menge $\til{\mathbf{L}}_d(R)$ aller $\F_q$-linearen Einbettungen 
$\theta\colon H\hookrightarrow R$ mit nicht-entartetem Bild zuordnet.
Dazu geh\"ort in der Kategorie der $k$-Schemata ein kontravarianter Funktor $\til{\mathcal{L}}_d$ gegeben durch die \'etale Garbifizierung von $$\mathbf{Sch}_k
\rightarrow\mathbf{Mengen} \ : \ X\mapsto \{ (\mathcal{V},e_1,\ldots,e_d) \, | \, \mathcal{V}\in\mathcal{L}_d(X) \}$$ wobei
$e_1,\ldots,e_d$ Schnitte sind, die eine $\F_q$-Basis bilden in jedem geometrischen Punkt, d.h.\ dem separabelen Abschluss von jedem Punkt auf $X$. 

\vv

\paragraph\label{subpara-ring1}
Sei $\mathfrak{R}_d=k[T_1,\ldots,T_d,T^{-1}_d]$ und betrachte das $\F_q$-lineare Polynom 
$$
P^{\ast}(T_1,\ldots,T_d)(X)=P(1,T_1,\ldots,T_d).
$$
Schreiben wir $\til{\mathfrak{R}}_d$ f\"ur den ganzen Abschluss von $\mathfrak{R}_d$ im Zerf\"allungsk\"orper von $P(T_1,\ldots,T_d)$ \"uber $k(T_1,\ldots,T_d)$, so ist  
 $\til{\mathfrak{R}}_d/\mathfrak{R}_d$ eine Galoiserweiterung mit Gruppe $\mathfrak{G}=\GL(d,\F_q)$, insbesondere ist $\til{\mathfrak{R}}^{\mathfrak{G}}_d=\mathfrak{R}_d$
(\ref{lem-covaff2}).
Proposition (\ref{lem-GLquot1}) impliziert, dass
$\til{\mathfrak{R}}_d$ gleich der lokalisierte Polynomring
$$
\til{\mathfrak{R}}_d\cong k[U_1,\ldots,U_d,W^{-1}]
$$
ist, wobei $U_1,\ldots,U_d$ linear unabh\"angige L\"osungen von $P^{\ast}(T_1,\ldots,T_d)$ sind und
$W$ das Produkt aller nicht-Null Elemente von $\til{V}=\sum_{i=1}^d\F_qU_i$ ist.
Man bemerke, dass $\til{V}$ offensichtlich zu 
$\mathbf{L}_d(\til{\mathfrak{R}}_d)$ geh\"ort.

\vv

\paragraph\label{subpara-functor2}
Offenbar wird der Funktor $\til{\mathbf{L}}_d$ vom Ring {\sl $\til{\mathfrak{R}}_d$ mit der Basis $U_1,\ldots,U_d$} dargestellt.
Was den Funktor $\mathbf{L}_d$ angeht, erwartet man wegen 
$\mathbf{L}_d=\mathfrak{G}\backslash\til{\mathbf{L}}_d$, dass 
$\mathfrak{R}_d$ der entsprechende Ring ist.
Sei $h_{\mathfrak{R}_d}=\Hom_k(\mathfrak{R}_d,\textrm{--})$ der von ihm dargestellte Funktor.
Es gibt einen nat\"urlichen Morphismus zwischen Funktoren
$
\Phi_d\colon \mathbf{L}_d\longrightarrow h_{\mathfrak{R}_d},
$
der jedem $V\in \mathbf{L}_d(R)$  die Abbildung $$\Phi_d(R)(V) \ \colon \  \mathfrak{R}_d\rightarrow R \ \colon \ T_i\mapsto a_i \ (i=1,\ldots,d)$$ zuordnet, wobei die $a_i$ dadurch definiert sind, dass $P^{\ast}(a_1,\ldots,a_d)$  das nach (\ref{lem-ring1}) zu $V$ geh\"orige $\F_q$-lineare Polynom 
\"uber $R$ ist. Dabei ist $\Phi_d$ ein Monomorphismus, aber im Allgemeinen kein Epimorphismus, da das Polynom $P^{\ast}(a_1,\ldots,a_d)$ nicht in $R$  zu zerfallen braucht --- deswegen stellt
$\mathfrak{R}_d$ den Funktor $\mathbf{L}_d$ nicht fein dar. Wir haben jedoch die

\vv

\paragraph\label{prop-functor2}
{\bf Proposition.}\ {\sl Das Paar $(\mathfrak{R}_d,\Phi_d)$ stellt den Funktor
$\mathbf{L}_d$ grob dar.}

\vv

\pf
F\"ur einen algebraisch abgeschlossenen K\"orper $F$ ist $\Phi_d(F)$ deutlich bijektiv.
Sei $S$ eine $k$-Algebra, und $\Psi\colon \mathbf{L}_d\rightarrow h_S$ ein Morphismus von Funktoren. Betrachte $\Psi(\til{\mathfrak{R}}_d)\colon \mathbf{L}_d(\til{\mathfrak{R}}_d)\rightarrow
\Hom_k(S,\til{\mathfrak{R}}_d)$, und sei $\phi$ das Bild von $\til{V}$.
Wegen Funktorialit\"at operiert die Galoisgruppe $\mathfrak{G}$ auf $\mathbf{L}_d(\til{\mathfrak{R}}_d)$ sowie $\Hom_k(S,\til{\mathfrak{R}}_d)$, und die Abbildung $\Psi(\til{\mathfrak{R}}_d)$ ist offensichtlich
$\mathfrak{G}$-\"aquivariant.
Da $\til{V}$ $\mathfrak{G}$-invariant ist als Element von
$\mathbf{L}_d(\til{\mathfrak{R}}_d)$, 
ist $\phi$ $\mathfrak{G}$-invariant, weswegen $\phi\in\Hom_k(S,\mathfrak{R}_d)$.
Es gibt also eine Abbildung $\phi^{\ast}\colon h_{\mathfrak{R}_d}\rightarrow h_S$.

Es muss jetzt gezeigt werden, dass $\phi^{\ast}\circ\Phi_d=\Psi$.
Sei $R$ eine $k$-Algebra, und betrachte die Abbildungen
$\phi^{\ast}(R)\circ\Phi_d(R)$ und $\Psi(R)$.
Sei $V\in \mathbf{L}_d(R)$ und $\psi\colon\mathfrak{R}_d\rightarrow R$ der zugeh\"orige Morphismus (d.h.\ $\psi=\Phi_d(R)(V)$). Definiere $\varphi=\Psi(R)(V)\colon S\rightarrow R$.
Man kann einen Morphismus  $\til{\mathfrak{R}}_d\rightarrow R$ w\"ahlen, der mit $\psi$
und $\mathfrak{R}_d\hookrightarrow\til{\mathfrak{R}}_d$ kompatibel ist und $\til{V}$ bijektiv auf $V$ abbildet. Dann kommutiert
$$
\begin{array}{ccc}
\mathbf{L}_d(\til{\mathfrak{R}}_d)&\longrightarrow&\Hom_k(S,\til{\mathfrak{R}}_d)\\
\bigdownarrow&&\bigdownarrow\\
\mathbf{L}_d(R)&\longrightarrow&\Hom_k(S,R)\rlap{,}
\end{array}
\qquad
\begin{array}{ccc}
\til{V}&\longmapsto&\phi\\
\setlength{\unitlength}{1pt}
\begin{picture}(9,22)(0,0)
\put(2,18){\line(1,0){4}}
\put(4,18){\vector(0,-1){16}}
\end{picture}&&
\setlength{\unitlength}{1pt}
\begin{picture}(9,22)(0,0)
\put(2,18){\line(1,0){4}}
\put(4,18){\vector(0,-1){16}}
\end{picture}\\
V&\longmapsto&\varphi\rlap{,}
\end{array}
$$
weswegen $\varphi=\psi\circ\phi$, und also folgt die erw\"unschte Gleichung $\phi^{\ast}(R)\circ\Phi_d(R)=\Psi(R)$.
\qed 

\vv

\paragraph\label{subpara-flagfunctor1} {\bf Variante: Moduln von Fahnen.}
\ Wir werden ab jetzt gelegentlich die Schreibweise $P \circ Q$ benutzen f\"ur die 
Komposition zweier Polynome: $P \circ Q (X) := P(Q(X))$. 

\vv

\paragraph\label{lem-flagfunctor1}
{\bf Lemma.}\ {\sl Seien $R$ eine $\F_q$-Algebra und $V$ und $V'$ nicht-entartete $\F_q$-lineare Unterr\"aume von $R$ der Dimension $d$, respektive $d'$, wobei wir $d \geq d'$ unterstellen. Seien $P$ und $P'$ die zu $V$ und $V'$ geh\"origen $\F_q$-linearen Polynome der respektiven H\"ohen $d,d'$. Dann sind folgende Eigenschaften \"aquivalent:

(i) Es gibt ein $\F_q$-lineares Polynom $P''$ der H\"ohe $d-d'$ mit $P = P'' \circ P'$. 

(ii) $V' \subseteq V$. }

\vv

\pf Da $V$ und $V'$ die respektiven Mengen der Wurzeln von $P$ und $P'$ sind wie in \ref{lem-ring1},  ist (i) $\Rightarrow $(ii) klar. Die andere Richting beweisen wir, indem wir ein $\F_q$-lineares Komplement $W \subseteq R$ von $V'$ in $V$ betrachten: $V=V' \oplus W$. Weil die Polynomabbildung $P' : W \rightarrow R$ eine $\F_q$-lineare Injektion ist, so ist ihr Bild $V''$ ein $\F_q$-linearer Unterraum von $R$ der Dimension $d-d'$. Sei $P'' = \prod_{v'' \in V''} (X-v'')$ das zu $V''$ geh\"orige Polynom. Wir k\"onnen voraussetzen, dass $P,P'$ und $P''$ monisch sind, und dann ist
\begin{eqnarray*}
P(X)&=&\prod_{w\in W}\prod_{v'\in V'}(X-v'-w)=\prod_{w\in W}P'(X-w)\\
&=&\prod_{w\in W}(P'(X)-P'(w))=P''\circ P'(X).
\end{eqnarray*}
Damit ist das Lemma bewiesen. \qed

\vv

\paragraph\label{cor-flagfunctor2}
{\bf Korollar.}\ {\sl Sei $d=f_1+\cdots+f_e$ eine Zerlegung  von $d$ in positive ganze Zahlen. F\"ur jede Fahne $V_1\subset V_2\subset\cdots\subset V_e$ bestehend aus nicht-entarteten $V_i$ der Dimension $d_i=f_1+\cdots+f_i$ gibt es ein bis auf Skalarmultiplikation mit $R^{\times}$ eindeutig bestimmtes nicht-entartetes $\F_q$-lineares Polynom \"uber $R$ der Form
$P_e\circ P_{e-1}\circ\cdots\circ P_1$, wobei $P_i$ nicht-entartet der H\"ohe  $f_i$ ist.}

\vv

\pf Per Induktion wird die Behauptung auf den Fall $e=2$, der im obigen Lemma behandelt wurde, zur\"uckgef\"uhrt. \qed

\vv

\paragraph\label{subpara-flagfunctor3} {\bf Parabolischer Quotientenfunktor.} \ F\"ur eine Zerlegung einer ganzen Zahl $d=f_1+\cdots+f_e$ in positive ganze Zahlen definieren wir einen kovarianten Funktor  $\mathbf{L}_{(f_1,\ldots,f_e)}\colon
\mathbf{Alg}_k\rightarrow\mathbf{Mengen}$ der jeder $k$-Algebra die Menge aller Fahnen $V_1\subset V_2\subset\cdots\subset V_e$ nicht-entarteter $\F_q$-linearer Unterr\"aume $V_i$ von $R$ der Dimension $d_i := f_1+\cdots+f_i$ zurordnet. Die zugeh\"orige \'etale Garbe $\mathcal{L}_{(f_1,\ldots,f_e)}$ ist gegeben durch $$
X\mapsto\mathcal{L}_{(f_1,\ldots,f_e)}(X)=
\left\{
\begin{minipage}{9em}
\setlength{\baselineskip}{.85\baselineskip}
\begin{small} \begin{center}
{\slshape Fahnen $\mathcal{V}_1\subset\cdots\subset\mathcal{V}_e$
mit $\mathcal{V}_i\in\mathcal{L}_{d_i}(X)$}
\end{center} \end{small}
\end{minipage}
\right\}.
$$
Offensichtlich ist der Funktor $\mathbf{L}_{(f_1,\ldots,f_e)}$ der Quotient von $\til{\mathbf{L}}_d$ durch die parabolische Untergruppe $\mathfrak{P}_{(f_1,\ldots,f_e)}$ von $\mathfrak{G}$ des richtigen Fahnentyps.

\vv

\paragraph\label{prop-flagfunctor3}
{\bf Proposition} (cf.\ \cite{Hew} 1.4).\ {\sl Mit der gleichen Notation wie (\ref{cor-GLquot1}) gilt 
\begin{eqnarray*}
& & \F_q[U_1,\ldots,U_d,W^{-1}]^{\mathfrak{P}_{(f_1,\ldots,f_e)}}= \\
& & \F_q[S^1_1,\ldots,S^1_{f_1},(S^1_{f_1})^{-1};
\ldots;S^e_1,\ldots,S^e_{f_e},(S^e_{f_e})^{-1}], \end{eqnarray*} wobei $P^{\ast}(S^e_1,\ldots,S^e_{f_e})\circ\cdots\circ P^{\ast}(S^1_1,\ldots,S^1_{f_1})=
X\prod_{v\in\til{V}-\{0\}}(1-v^{-1}X)$.}

\vv

\pf Wie (\ref{lem-GLquot1}). \qed

\vv

\paragraph \label{subpara-flagfunctor4} {\bf Darstellbarkeit.} \ Wir definieren $$\mathfrak{R}_{(f_1,\ldots,f_e)}=k[S^1_1,\ldots,S^1_{f_1},(S^1_{f_1})^{-1};
\ldots;S^e_1,\ldots,S^e_{f_e},(S^e_{f_e})^{-1}].$$
Nach \ref{cor-flagfunctor2} gibt es wie in (\ref{subpara-functor2}) einen nat\"urlichen Morphismus zwischen Funktoren 
$$
\Phi_{(f_1,\ldots,f_e)}\colon \mathbf{L}_{(f_1,\ldots,f_e)}\longrightarrow
h_{\mathfrak{R}_{(f_1,\ldots,f_e)}}.
$$

\vv

\paragraph\label{prop-flagfunctor4}
{\bf Proposition.}\ {\sl Das Paar $(\mathfrak{R}_{(f_1,\ldots,f_e)},\Phi_{(f_1,\ldots,f_e)})$
stellt den Funktor $\mathbf{L}_{(f_1,\ldots,f_e)}$ grob da.}

\vv

\pf Wie (\ref{prop-functor2}). \qed 

\vv

\paragraph\label{subpara-functorrelation2} {\bf Zusammenhang zwischen $\mathfrak{R}_d$ und $\mathfrak{R}_{(1,\ldots,1)}$.} \ 
Wir berechnen die Abbildung
$$
\mathfrak{R}_d=k[T_1,\ldots,T_d,T^{-1}_d]\longhookrightarrow
\mathfrak{R}_{(d-1,1)}=k[T'_1,\ldots,T'_{d-1},T^{\prime -1}_{d-1};S,S^{-1}]
$$
definiert durch $P^{\ast}(S)\circ P^{\ast}(T'_1,\ldots,T'_{d-1})=P^{\ast}(T_1,\ldots,T_d)$.
Weil
\begin{eqnarray*}
P^{\ast}(S)\circ P^{\ast}(T'_1,&\ldots&,T'_{d-1})=P(1,S)\circ P(1,T',\ldots,T'_{d-1})\\
&=&(SX^q+X)\circ P(1,T',\ldots,T'_{d-1})\\
&=&S\cdot P(1,T',\ldots,T'_{d-1})^q+P(1,T',\ldots,T'_{d-1})\\
&=&P^{\ast}(S+T'_1,ST^{\prime q}_1+T'_2,\ldots,ST^{\prime q}_{d-2}+T'_{d-1},ST^{\prime q}_{d-1}),
\end{eqnarray*}
findet man $\mathfrak{R}_{(d-1,1)}=\mathfrak{R}_d[S]$ und das Minimalpolynom von $S$
ist gegeben durch das monische Polynom $$
S^{\frac{q^d-1}{q-1}}+\sum^{d-1}_{i=0}(-1)^{d-i}T^{q^i}_{d-i}S^{\frac{q^i-1}{q-1}}.
$$
Per Induktion folgt hieraus, dass die Erweiterung 
$$
\mathfrak{R}_d=k[T_1,\ldots,T_d,T^{-1}_d]\longhookrightarrow
\mathfrak{R}_{(1,\ldots,1)}=k[S_1^{\pm 1},\ldots,S_d^{\pm 1}]
$$
gegeben durch $P^{\ast}(S_d)\circ\cdots\circ P^{\ast}(S_1)=P^{\ast}(T_1,\ldots,T_d)$  endlich von Grad $\prod^d_{i=1}\frac{q^i-1}{q-1}$ ist.

\vv

\paragraph\label{subpara-functorelation3} {\bf Zusammenhang zwischen $\mathfrak{R}_{(1,\ldots,1)}$ und $\til{\mathfrak{R}}_d$.} 
\ Sei $\mathfrak{R}_{(1,\ldots,1)}=k[S_1^{\pm 1},\ldots,S_d^{\pm 1}]$, wobei die Parameter korrespondieren zu denen in $$P^{\ast}(S_d)\circ\cdots\circ P^{\ast}(S_1),$$ und sei $U_1,\ldots,U_d$ die Basis aus  (\ref{subpara-ring1}), derart dass der von $U_1,\ldots,U_i$ \"uber $\F_q$ erzeugte Raum 
die Menge der Wurzeln von $$P^{\ast}(S_i)\circ\cdots\circ P^{\ast}(S_1)$$
ist f\"ur $i=1,\ldots,d$.
Dann gilt induktiv 
$$
\begin{array}{c}
U_1^{q-1}+S^{-1}_1=0,\\
\{P^{\ast}(S_1)(U_2)\}^{q-1}+S^{-1}_2=0,\\
\vdots\\
\{P^{\ast}(S_{d-1})\circ\cdots\circ P^{\ast}(S_1)(U_d)\}^{q-1}+S^{-1}_d=0,
\end{array}
$$
welches die definierende Gleichungen sind f\"ur den Morphismus $$\mathfrak{R}_{(1,\ldots,1)}\hookrightarrow\til{\mathfrak{R}}_d$$ der also endlich von Grad $q^{\frac{d(d-1)}{2}}(q-1)^d$ ist.

\vv

\vv

\sectioning{Konstruktion einer Familie gez\"ugelter \"Uberlagerungen}\label{familie}

\vv

\paragraph\label{prop-invariant}
{\bf Proposition} (Konstruktion gez\"ugelter Erweiterungen). \ {\sl (i) Sei $A$ komplett, $B/A$ eine gez\"ugelte Erweiterung vom Typ $(n,d)$ mit Galoisgruppe $G$ und $H$ die $p$-Sylowgruppe von $G$. W\"ahle eine Uniformisierende $x$ f\"ur $B$, derart dass die Operation von $G$ auf $B$ gegeben ist in der kanonischen Form (\ref{prop-canonicalform}) und sei $\theta=\theta_{B/A}/dx$. 
Sei $P^{\ast}(a_1,\ldots,a_d)$ ($a_d\neq 0$) das $\F_q$-lineare Polynom, das zu $V=\theta(H)$ geh\"ort wie in (\ref{lem-ring1} mit $R=k$). Betrachte das Element $y$ in dem Quotientenk\"orper $L$ von $B$, das durch $$P^{\ast}(a_1,\ldots,a_d)(x^{-1})^n=y^{-1}$$ definiert wird, d.h.\ $Q(x) = 0$ mit
$$
Q(x)=x^{nq^d}-y(x^{q^d-1}+a_1x^{q(q^{d-1}-1)}+\cdots+a_{d-1}x^{q^{d-1}(q-1)}+a_d)^n.
$$
Dann ist $y \in A$ eine Uniformisierende, und $B=A[x]/(Q(x))$.

(ii) Sei umgekehrt $(a_1,\ldots,a_d) \in k^d$ mit $a_d \neq 0$ und definiere $B=A[x]/(Q(x))$, wobei das Polynom $Q$ wie oben gegeben ist. Dann ist $B/A$ eine gez\"ugelte Erweiterung vom Typ $(n,d)$ mit Galoisgruppe $G$, und $x$ eine Uniformisierende f\"ur $B$, derart dass $(\theta_{B/A}/dx)(H)$ (mit $H$ die $p$-Sylowgruppe von $G$) der $\F_p$-lineare Unterraum von $k$ ist, der durch die Gleichung $$P^*(a_1,\ldots,a_d)(X)=0$$ bestimmt wird. }

\vv

\pf (i) Wegen $Q(x)=0$ und $a_d \neq 0$ ist $y \in B$. Da die Operation von $G$ durch $\gamma(y^{-1})=\zeta^{-1}y^{-1}$ und $\sigma(y^{-1})=y^{-1}+\theta(\sigma)$ f\"ur $\sigma\in H$ gegeben wird, sehen wir, dass $y$ $G$-invariant ist, also $y \in B^G=A$. Da $B/A$ totalverzweigt ist vom Grad $nq^d$ sind Uniformisierenden von $A$ Elemente aus $A$ mit $B$-Bewertung $nq^d$; dies ist der Fall f\"ur $y$. Sei $B'=A[x]/(Q(x))$. Nach (\cite{Ser} I, Prop.\ 17) ist $B'$ ein diskreter Bewertungsring mit Uniformisierender $x$ und Restklassenk\"orper $k$. Weil $B' \subseteq B$ und die Quotientenk\"orper von $B$ und $B'$ denselben Grad haben, also gleich sind, gilt auch $B'=B$. 

(ii) Weil $Q(x)$ Eisensteinsch ist, ist $B$ ein diskreter Bewertungsring mit $x$ als Uniformisierender, und $B/A$ ist dabei totalverzweigt (\cite{Ser} I, Prop.\ 17). Da $Q(x)$ invariant ist unter den Operationen  $y\mapsto\zeta y$ und  
$y\mapsto\frac{y}{1-uy}$, wobei $u$ den durch $P^*(a_1,\ldots,a_d)(X)=0$ bestimmten $\F_q$-linearen Unterraum von $k$ durchl\"auft, ist $B/A$ gez\"ugelt vom Typ $(n,d)$. \qed

\vv

\paragraph\label{cor-isomaction1}
{\bf Korollar.}\ {\sl F\"ur komplette $A$ gibt es eine nat\"urliche Bijektion $$
\left\{
\begin{minipage}{11em}
\setlength{\baselineskip}{.85\baselineskip}
\begin{small} \begin{center}
{\slshape Isomorphieklassen gez\"ugelter Erweiterungen $B/A$ vom Typ $(n,d)$}
\end{center} \end{small} 
\end{minipage}
\right\}
\stackrel{\sim}{\longrightarrow}
\GL(d,\F_q)\backslash\left\{
\begin{minipage}{7.5em}
\setlength{\baselineskip}{.85\baselineskip}
\begin{small}
{\slshape $\F_q$-lineare $H\hookrightarrow k$}
\end{small}
\end{minipage}
\right\}/k^{\times}.
$$}

\pf Direkt aus (\ref{prop-isomaction}) und  (\ref{prop-invariant})
(2). \qed

\vv

\paragraph\label{cor-isomaction2}
{\bf Korollar.}\ {\sl F\"ur $d=0$ und $1$ gibt es nur eine gez\"ugelte Erweiterung von Typ $(n,d)$ eines gegebenen kompletten Ringes $A$.}

\vv

\pf Folgt aus (\ref{cor-isomaction1}).  
F\"ur $d=0$ ist das Ergebniss aber auch klassisch, f\"ur $d=1$ folgt es 
unmittelbar aus der Artin-Schreierschen Theorie. \qed

\vv

\paragraph\label{subpara-genext1} {\bf Konstruktion des Schemas $\mathcal{S}$.} \ Sei $k$ einen K\"orper \"uber $\F_q$ und $A$ ein diskreter Bewertungsring mit Restklassenk\"orper $k$. Nach (\ref{cor-isomaction1}) muss die allgemeinste Familie gez\"ugelter Erweiterungen von $A$ \"uber $$\mathfrak{G}_k\backslash\til{\mathbf{L}}_d/\G_{m,k}=
\mathbf{L}_d/\G_{m,k}$$ konstruiert werden, wo $\mathfrak{G}=\GL(d,\F_q)$. Der Funktor $\mathbf{L}_d$ wird grob vom Ring $\mathfrak{R}_d=
k[T_1,\ldots,T_d,T^{-1}_d]$ dargestellt (\ref{prop-functor2}), wobei die Parameter aus dem Polynom $P^{\ast}(T_1,\ldots,T_d)=P(1,T_1.\ldots,T_d)$ stammen.
Sei $$\til{\mathfrak{R}}_d=k[U_1,\ldots,U_d,W^{-1}]$$ wie in (\ref{subpara-ring1}) und $\til{V}=\sum^d_{i=1}\F_qU_i$ der Vektorraum der Wurzeln von  
$$P^{\ast}(T_1,\ldots,T_d).$$ 
F\"ur $c\in\til{\mathfrak{R}}^{\times}_d$ hat $c^{-1}\til{V}$ als zugeh\"origes Polynom  
$P^{\ast}(c^{q-1}T_1,\ldots,c^{q^d-1}T_d)$, weil $$
P^{\ast}(T_1,\ldots,T_d)(cX)=cP^{\ast}(c^{q-1}T_1,\ldots,c^{q^d-1}T_d)(X).
$$
Wir versehen also $\mathfrak{R}_d$ mit einer Graduierung, deren Gewichte $w(T_i)=q^i-1$ ($i=1,\ldots,
d$) sind. Dann ergibt sich sofort, dass der Funktor $\mathbf{L}_d/\G_{m,k}$ grob dargestellt wird von (dem Koordinatenring von) dem affinen offenen Unterschema $$\mathcal{S} \subset \P_k(q-1,\ldots,q^d-1)=\Proj 
k[T_1,\ldots,T_d],$$ welches das Komplement der Hyperebene $V(T_d)$ ist.

\vv

\paragraph\label{subpara-genext2} {\bf Konstruktion des Schemas $\mathcal{X}$.} \ 
Man betrachte jetzt den graduierten Ring $k[T_1,\ldots,T_d,X,Z]$, mit Gewichten $w(X)=w(Z)=1$ und, f\"ur eine feste Uniformisierende $y$ in $A$ (cf.\ \ref{prop-invariant}), das Polynom
$$
Q=X^{nq^d}-yZ^n(X^{q^d-1}+T_1X^{q(q^{d-1}-1)}+\cdots+T_{d-1}X^{q^{d-1}(q-1)}+T_d)^n,
$$
welches homogen ist vom Grad $nq^d$ in $A[T_1,\ldots,T_d,X,Z]$.
Es definiert ein geschlossenes Teilschema $V(Q)$ in 
$$\P_A(q-1,\ldots,q^d-1,1,1)=\Proj A[T_0,\ldots,T_d,X,Z].$$
Sei $\mathcal{X}$ das Komplement von $V(T_d)$ in $V(Q)$.
Wir erhalten dann einen Morphismus von $A$-schemata
$$
\pi\colon\mathcal{X}\longrightarrow\mathcal{S}\times_k\Spec A,
$$
der endlich vom Grad $nq^d$ ist.
F\"ur komplettes $A$ und einen algebraisch abgeschlossenen K\"orper $F/k$ entsteht jede gez\"ugelte Erweiterung von $A\otimes_kF$ durch geschickte Spezialisierung von $T_1,\ldots, T_d$ und $Z$ \"uber $k$, wobei das Bild von $T_d$ und $Z$ nicht-null ist.

\vv

\paragraph\label{subpara-genext3}  {\bf Konstruktion des Schemas $\til{\mathcal{S}}$.} \  Zwischen $\til{\mathbf{L}}_d$ und $\mathbf{L}_d$ liegt der "`vollst\"andige-Fahnen"'-Funktor $\mathbf{L}_{(1,\ldots,1)}$, der vom Ring
$$\mathfrak{R}_{(1,\ldots,1)}=k[S_1^{\pm 1},\ldots S_d^{\pm 1}]$$ grob dargestellt wird, wobei $S_1,\ldots,S_d$ mit dem Parameter im Polynom $P^{\ast}(S_d)\circ\cdots\circ P^{\ast}(S_1)$ \"ubereinstimmen.
Weil $P^{\ast}(S)(cX)=cP^{\ast}(c^{q-1}S)(X)$, gilt $$
P^{\ast}(S_d)\circ\cdots\circ P^{\ast}(S_1)(cX)=
cP^{\ast}(c^{q-1}S_d)\circ\cdots\circ P^{\ast}(c^{q-1}S_1)(X),
$$
und deswegen ist $w(S_i)=q-1$.
\"Ahnlich wie oben wird der Funktor $\mathbf{L}_{(1,\ldots,1)}/\G_{m,k}$ also grob dargestellt durch das affine Unterschema von  $$\P_k(q-1,\ldots,q-1)=\Proj k[S_1,\ldots,S_d]$$ welches das Komplement der Hyperebenen $\bigcup_{i=1}^d V(S_i)$ ist.
Aus (\ref{subpara-functorelation3}) folgt au{\ss}erdem, dass im Ring $\til{\mathfrak{R}}_d
=k[U_1,\ldots,U_d,W^{-1}]$ die Gewichte  $w(U_i)=-1$ sind f\"ur $i=1,\ldots,d$.
Deswegen wird $\til{\mathbf{L}}_d/\G_{m,k}$ grob dargestellt durch das affine offene Unterschema von $$\P^{d-1}_k=\Proj k[U^{-1}_1,\ldots,U^{-1}_d]$$ welches das Komplement der Hyperebenen $\bigcup_{i=1}^d V(U^{-1}_i)$ und aller Hyperfl\"achen definiert durch nicht-triviale $\F_q$-lineare Kombinationen von $$U_1, \ldots, U_d$$ ist. Sei $\til{\mathcal{S}}$ das derart konstruierte $k$-Schema. 
\vv

\paragraph\label{subpara-genext4} {\bf Konstruktion des Schemas $\til{\mathcal{X}}$.} \ 
Es gibt einen endlichen Morphismus $\til{\mathcal{S}}\rightarrow\mathcal{S}$, der Galois ist mit Galoisgruppe $\mathfrak{G}=\GL(d,\F_q)$.
Es sei $\til{\pi}\colon\til{\mathcal{X}}\rightarrow\til{\mathcal{S}}$ durch Kommutieren des Diagrams $$
\begin{array}{ccc}
\til{\mathcal{X}}&\longrightarrow&\mathcal{X}\\
\llap{$\scriptstyle{\til{\pi}}$}\bigdownarrow&&\bigdownarrow\rlap{$\scriptstyle{\pi}$}\\
\til{\mathcal{S}}\times_k\Spec A&\longrightarrow&
\mathcal{S}\times_k\Spec A\rlap{.}
\end{array}
$$
definiert. 
Dieser Morphismus wird durch Zur\"uckziehen von $Q$ bekommen, ist also gegeben durch
$$yZ^nX^{-n}\prod_{v\in\til{V}-\{0\}}(1-v^{-1}X^{-1})^n=1,$$ was gleichbedeutend damit ist, dass $$\til{Q}=0 \mbox{  \ \ wo  \ \ } \til{Q}=X^{nq^d}-yZ^n\prod_{v\in\til{V}-\{0\}}(X-v^{-1})^n
$$
(man bemerke dass jedes $v\in\til{V}-\{0\}$ Gewicht $-1$ hat).
So erlaubt $\til{\mathcal{X}}$ eine Operation durch $G_{n,d}$ gegeben von 
$$
 \gamma(X\colon Z)=(X\colon\zeta^{-1}Z),\qquad
\sigma(X\colon Z)=(X: Z-\theta(\sigma)X), $$
wobei $\theta\colon H\stackrel{\sim}{\rightarrow}\til{V}$ irgendein $\F_q$-linearer Isomorphismus und $\gamma$ erzeugendes Element eines Torusteils ist. 
Diese Definition ist sinnvoll, weil $\til{Q}$ nicht Null ist f\"ur $X=v^{-1}$ ($v\in\til{V}-\{0\}$) und man also  $1-\theta(\sigma)Y$ invertieren kann auf $\til{\mathcal{X}}$. Offensichtlich ist dann $\til{\pi}$ eine Galoissche \"Uberlagerung mit der Gruppe $G_{n,d}$, und kann jede gez\"ugelte Erweiterung \"uber $A\otimes_kF$ (mit $F$ einen algebraisch abgeschlossenen K\"orper) durch Spezialisierung \"uber $k$ bekommen werden, und zwar eindeutig bis auf $\mathfrak{G}$.

\vv

\paragraph\label{subpara-deg1} \ {\bf Der Funktor $\mathbf{L}_d^1$.} \ Damit wir jetzt eine Entartung konstruieren k\"onnen, nehmen wir an, dass $d>1$, und betrachten den Funktor $\mathbf{L}^1_d$ der nicht-entartetem $\F_q$-linearem Unterr\"aume der Dimension $d$ {\sl mit einem ausgezeichneten nicht-null Vektor} klassifiziert. Der Funktor $\mathbf{L}^1_d$ ist eine \"Uberlagerung vom Grad $q-1$ \"uber $\mathbf{L}_{(1,d-1)}$ (siehe oben). Der Ring $$\mathfrak{R}_{(1,d-1)}=k[S,S^{-1},S_1,\ldots,S_{d-1},S^{-1}_{d-1}]$$
stellt $\mathbf{L}_{(1,d-1)}$ grob da (\ref{subpara-flagfunctor4}), ge\-h\"o\-rend zum Polynom 
$$
P^{\ast}(S_1,\ldots,S_{d-1})\circ P^{\ast}(S).
$$
Deswegen wird der Funktor $\mathbf{L}^1_d$ vom Ring $$\mathfrak{R}^1_d=k[T,T^{-1},S_1,\ldots,S_{d-1},S^{-1}_{d-1}],$$ grob dargestellt, wobei $T$ durch $$
T^{q-1}+S^{-1}=0
$$
bestimmt wird (cf.\ \ref{subpara-functorelation3}).
 
Aus (\ref{prop-invariant}) folgt, dass eine Familie gez\"ugelter Erweiterungen \"uber 
$\mathfrak{R}^1_d$ gegeben wird durch die Gleichung
$$
(P^{\ast}(S_1,\ldots,S_{d-1})(T^{1-q}x^{-q}+x^{-1}))^n=y^{-1},
$$ oder auch durch ein Eisensteinsches Polynom wie in
(\ref{prop-invariant}).
Wir untersuchen was passiert, wenn $T\rightarrow\infty$, d.h.\ wenn der ausgezeichnete Vektor $v\rightarrow 0$, w\"ahrend die \"ubrigen Parameter $$(S_1,\ldots,S_{d-1})=(s_1,\ldots,s_{d-1})$$ f\"ur {\sl feste} $s_i\in k$ erf\"ullen. Es sei ab jetzt $k$ algebraisch abgeschlossen.
Wir setzen $t=-T^{-1}$ ein, und ersetzen $y$ durch $-y$.
Wir betrachten die von $\G_{m,k}=\Spec k[t^{\pm 1}]$ parametrisierte Familie $\mathcal{X}_t$ gez\"ugelter Erweiterungen definiert durch 
 $$P^{\ast}(s_1,\ldots,s_{d-1})(t^{q-1}x^{-q}-x^{-1})^n=
y^{-1}$$ in $\Spec A[x,t]$.
Die Operation der Galoisgruppe $G_{n,d}$ sieht dann so aus:
man betrachtet eine $\F_q$-lineare Zerlegung $H=H'\oplus H''$ der $p$-Sylowgruppe $H$ von $G_{n,d}$ und w\"ahlt die Einbettung $\theta\colon H\rightarrow k[t^{\pm 1}]$ so, dass 
$\theta(H)$ zu $P^{\ast}(s_1,\ldots,s_{d-1})\circ P^{\ast}(-t^{q-1})$ geh\"ort und
$\theta(H'')$ \"uber $\F_q$ von $t$ erzeugt wird.
Sei $\gamma$ ein erzeugendes Element eines Torusteils mit $\gamma\sigma\gamma^{-1}=
\zeta\sigma$ f\"ur $\sigma\in H$.
Dann wird die Galoissche Operation wie in (\ref{prop-canonicalform}) gegeben.

\vv

\paragraph\label{prop-deg-vorbereitung} \ {\bf Definition/Konstruktion.} \ Sei $\mathfrak{R}$ der ganze Abschluss von $k[t^{\pm 1}]$ im Zerf\"allungsk\"orper des Polynoms $$P^{\ast}(t^{1-q}s_1,\ldots,t^{1-q^d}s_{d-1})\circ P^{\ast}(-1)$$ \"uber 
$k(t)$. F\"ur jede L\"osung $u$ von $P^{\ast}(s_1,\ldots,s_{d-1})$ geh\"oren dann alle Wurzeln $w$ von
$w^q-w-tu$ zu $\mathfrak{R}$.
Sei $V$ der $\F_q$-lineare Unterraum von 
$\mathfrak{R}$ gegeben durch alle solche Wurzeln $w$, und bezeichne $V''=\F_q\subseteq V$.
Sei $\eta\ \colon \ V\rightarrow k$ die Abbildung $w\mapsto t^{-1}(w^q-w)$. Dann hat $\eta$ als Bild genau den zu $P^{\ast}(s_1,\ldots,s_{d-1})$ geh\"orenden $\F_q$-linearen Unterraum. Man nimmt an\-schlie\-{\ss}end einen Schnitt von $\eta$, und es ergibt sich eine Zerlegung $V=V''\oplus V'$, derart dass $\eta$   den Raum $V'$ bijektiv auf die Wurzeln von $P^{\ast}(s_1,\ldots,s_{d-1})$ abbildet. Sei $\til{\theta}=\theta'\oplus\theta''\colon H=H'\oplus H''\rightarrow V=V'\oplus V''$ eine 
$\F_q$-lineare Einbettung die diese Zerlegung respektiert.

\vv

\paragraph\label{prop-deg3}
{\bf Proposition.}\ {\sl Wir setzen $d>1$ voraus. Dann gilt folgendes. 

(i) Die Familie $\mathcal{X}_t$ ist \"uber $\Spec A\times_k\G_{m,k}=\Spec A[t^{\pm 1}]$ isomorph zur Familie $\til{\mathcal{X}}_t$, die in  
$\Spec A[y,t]\times_A\P^1_A=\Spec A[y,t]\times_A\Proj A[Z,W]$  durch folgende Gleichungen definiert wird $$
\left\{
\begin{array}{l}
P^{\ast}(s_1,\ldots,s_{d-1})(x^{-1})^n=y^{-1},\\
x(Z^{-q}-Z^{-1}W^{1-q})=tW^{-q}.
\end{array}
\right.
$$

(ii) Der Isomorphismus zwischen $\mathcal{X}_t$ und $\til{\mathcal{X}}_t$ kann \"uber eine geschickte Erweiterung von $k[t^{\pm 1}]$ zu einem Isomorphismus von $G_{n,d}$-Galois\-er\-wei\-terung\-en \"uber $\Spec A[t^{\pm 1}]$ gemacht werden:
mehr spezifisch ist die Operation von $G_{n,d}$ auf $\til{\mathcal{X}}_{t,\mathfrak{R}}$ definiert durch
$$
\left\{
\begin{array}{l}
\gamma(x)=\zeta x,\\
\gamma(Z\colon W)=(Z\colon \zeta^{-1}W),
\end{array}
\right.
\qquad
\left\{
\begin{array}{l}
\sigma(x)={\textstyle \displaystyle \frac{x}{1-\eta\circ\theta'(\sigma')x}},\\
\sigma(Z\colon W)=(Z,Z+\til{\theta}(\sigma)W),
\end{array}
\right.
$$
wobei $\sigma=\sigma'+\sigma''\in H=H'\oplus H''$.

(iii) Die Faser $\til{\mathcal{X}}_0$ besteht aus einem $\P^1_k$ definiert durch $x=0$ (und in Folge dessen $y=0$) und $q$ Kopien von $\Spec B$, wobei $B/A$ die durch $$(P^{\ast}(s_1,\ldots,s_{d-1})(x^{-1}))^n=y^{-1}$$ definierte gez\"ugelte Erweiterung ist,  die $\P^1_k$ in den $q$ Punkten $\{ W/Z\in\F_q\}$ normal schneiden.

(iv) Die Gruppenoperation setzt sich auf nat\"urliche Weise auf der Faser $\til{\mathcal{X}}_0$ fort: $H'$ operiert wie \"ublich auf jedem $B$ und permutiert die $q$ Punkte $\{ W/Z\in\F_q\}$ und die zugeh\"orige Kopien von $\Spec B$.
Die Komponente $\P^1_k$ wird punktweise von $H'\rtimes T=G_{n,d-1}$ festgehalten, und der Punkt $\{ W/Z=\infty\}$ ist der einzige Punkt, der von ganz $G_{n,d}$ festgehalten wird.}

\vv

\pf (i) Es sei $I$ das Ideal 
$$
I=A[x,t]\cap(yP^{\ast}(s_1,\ldots,s_{d-1})(t^{q-1}x^{-q}-x^{-1})^n-1),
$$
wobei 
$(yP^{\ast}(s_1,\ldots,s_{d-1})(t^{q-1}x^{-q}-x^{-1})^n-1)$ Ideal ist in $A[x^{\pm 1},t],$ und es sei $J$ das Ideal
$$
J=A[x,t][Z/W]\cap(yP^{\ast}(s_1,\ldots,s_{d-1})(x^{-1})^n-1,x((Z/W)^{-q}-(Z/W)^{-1})-t),
$$
wobei der Schnitt in $A[x^{\pm 1},t][(Z/W)^{\pm 1}]$ genommen wird.
Es wird dann eine rationale Abbildung $\til{\mathcal{X}}_t\dashrightarrow\mathcal{X}_t$ definiert durch $$
\varphi\colon A[x,t]/I\longrightarrow A[x,t][Z/W]/J
$$
mit $\varphi(x)=t(Z/W)$. Wenn $t$ invertierbar ist, so besitzt $\varphi$ eine Umkehrung definiert durch $$\varphi^{-1}(x^{-1})=t^qx^{-q}-x^{-1} \ , \ \varphi^{-1}(Z/W)=t^{-1}x. $$

(ii) Dieser Morphismus ist $G_{n,d}$-\"aquivariant. 

(iii) Offensichtlich.

(iv) In der Tat bildet $\eta\circ\theta'$ den Raum $H'$ konstant auf den $\F_q$-linearen Raum ab, der zu $P^{\ast}(s_1,\ldots,s_{d-1})$ geh\"ort. Weil $V|_{t=0}=V\cap k=V''=\F_q$, gilt $\sigma(Z\colon W)=(Z,Z+\theta''(\sigma'')W)$ f\"ur $t=0$. \qed

\vv

\vv

\sectioning{Beweis f\"ur Satz A.}\label{beweis}

\vv

\paragraph \label{beweis-setup} {\bf Notationen.} \  Es sei jetzt $\pi : X \rightarrow Y$ eine gez\"ugelte Galoissche \"Uberlagerung glatter projektiver Kurven \"uber einen algebraisch abgeschlossenen K\"orper $k$ der Charakteristik $p>0$ mit Galoisgruppe $G$, verzweigt \"uber eine Menge endlich vieler abgeschlossener Punkte $B =\{ y_0, \ldots, y_r \} \hookrightarrow Y$, mit $y_0=y=\pi(x)$ f\"ur $x \in X$. Wir nehmen desweiteren an, $x$ habe eine Zerlegungsgruppe $G_x$ vom Typ $(n,d)$ mit $n=n_x$ und $q^d=p^{t_x}$. Man bemerke, dass $s_x \neq t_x$ gleichbedeutend damit ist, dass $d>1$.

\vv

\paragraph \label{local-global} {\bf Lokale und globale Funktoren.} \ Es sei 
$$\mathcal{F}_{(Y,G,(y_i,G_i)_{i=0}^r)} : \mathbf{Sch}_k \rightarrow \mathbf{Mengen}$$ der "`globale'' Funktor, der ein Schema $S$ \"uber $k$ abbildet auf die Menge gez\"ugelter \"Uberlagerungen von glatten projektiven Kurven \"uber $S$:
\unitlength0.7pt
\begin{center}
\begin{picture}(100,45)
\put(15,35){$\Pi$}
\put(-10,25){$\X$}
\put(4,30){\vector(1,0){34}}
\put(40,25){$\mathcal{Y} = S \times Y$}
\put(0,20){\vector(1,-1){15}}
\put(42,20){\vector(-1,-1){15}}
\put(43,6){$\scriptstyle{\mathrm{pr}}_1$}
\put(130,10){,} 
\put(15,-10){$S$}
\end{picture}
\end{center}
wobei es einen festen Punkt $s \in S$ gibt, derart dass $\Pi$ genau \"uber $\{ s \} \times B$ (gez\"ugelt) verzweigt und die Verzweigungsgruppe in $(s,y_i)$ isomorph zu $G_i$ ist. 

Wir ben\"otigen auch den "`lokalen"' Funktor $\mathcal{F}_{G} : \mathbf{Sch}_k \rightarrow \mathbf{Mengen}$ der einem $k$-Schema $S$ die Menge aller Galoisschen gez\"ugelten \"Uberlagerungen $T \rightarrow \Spec k [[y]] \times_k S$  mit Galoisgruppe $G$ zuordnet. 

Wie in \cite{BM}, 3.3. gibt es einen Morphismus von Funktoren 
$$ \mathcal{F}_{(Y,G,(y_i,G_i)_{i=0}^r)} \stackrel{\Psi}{\rightarrow} \prod_{i=0}^r \mathcal{F}_{G_i}.$$

\vv

\paragraph \label{formal-funktor} {\bf Formale Komplettierung.} \ Durch formale Komplettierung in einem $F$-wertigen Punkt (f\"ur beliebige K\"orper $F$ \"uber $k$), die wir mit 
\unitlength0.7pt
\begin{center}
\begin{picture}(100,45)
\put(15,35){$\pi'$}
\put(-15,25){$X'$}
\put(4,30){\vector(1,0){34}}
\put(43,25){$Y'=Y_F$}
\put(0,20){\vector(1,-1){15}}
\put(42,20){\vector(-1,-1){15}}
\put(1,-10){$\Spec F$}
\end{picture}
\end{center}
bezeichnen,  gehen die Funktoren und der Morphismus $\Psi$ in ihre formalen
Varianten \"uber, d.h.\ werden Deformationsfunktoren im Sinne {\sc Schlessinger}s 
$$ \hat{\mathcal{F}}_{(Y,G,(y_i,G_i)_{i=0}^r)} \stackrel{\hat{\Psi}}{\rightarrow} \prod_{i=0}^r \hat{\mathcal{F}}_{G_i}$$
in der Kategorie $\mathcal{C}_k$ der lokalen Artinschen $k$-Algebren mit Rest\-klassen\-k\"or\-per $k$. Dieser Morphismus wurde in \cite{BM} und \cite{CK2} betrachtet. 

Genauer gesagt haben wir Folgendes (wobei wir f\"ur Definitionen auf den ersten Abschnitt von \cite{CK2} verweisen): Sei $\rho': G \rightarrow \Aut_{F}(X')$ der Homomorphismus, der die Operation von $G$ auf $X'$ beschreibt. Sei $D_{X',\rho'}$ der Funktor aus \cite{CK2} 1.1, der {\sl alle} \"aquivarianten Deformationen  von $(X',\rho')$ zu $\mathcal{C}_k$ beschreibt. Dann ist der formale Funktor $\hat{\mathcal{F}}_{(Y,G,(y_i,G_i)^r_{i=0})}$ ein Teilfunktor von $D_{X',\rho'}$, weswegen der Tangentialraum von $\hat{\mathcal{F}}_{(Y,G,(y_i,G_i)^r_{i=0})}$ ein Unterraum von $D_{X',\rho'}(F[\varepsilon])$ ist.  Es seien $D_{\rho_i}$ die lokalen Funktoren aus \cite{CK2}, 1.11, die Deformationen zu $\mathcal{C}_k$ von der Operation der Verzweigungsgruppen auf den komplettierten lokalen Ringe $\hat{\mathcal{O}}_{X',x}$ beschreiben. Wie in \cite{CK2}, 1.6-1.9, wird der Zusammenhang zwischen Deformationen von $Y_F$, lokalen Deformationen der Verzweigungspunkte und Deformationen von $(X',\rho')$ beschrieben durch eine exakte Sequenz $$
0\longrightarrow\mathrm{H}^1(Y_F,\mathcal{T}_{Y}(-\Delta))\longrightarrow D_{X',\rho'}(F[\varepsilon])\longrightarrow\bigoplus^r_{i=0}D_{\rho'_i}(F[\varepsilon])\longrightarrow 0,
$$
wobei $\Delta$ der Divisor $$\Delta = \sum_{y \in T} y + \sum_{y \in W} 2y$$ ist, mit $T \cup W$ die Menge der Punkte, \"uber denen $\pi'$ verzweigt, und $y \in T$, genau dann, wenn $t_y = 0$ oder [$p=2$ und $t_y=1$]. F\"ur jede Deformation $(\til{X}',\til{\rho}')$ von $(X',\rho')$ liefert $\til{X}'/\til{\rho}'(G)$ eine Deformation von $Y_F$, weswegen es einen
Retrakt $$
\lambda\colon      D_{X',\rho'}(F[\varepsilon])\longrightarrow\mathrm{H}^1(Y_F,\mathcal{T}_{Y}(-\Delta)),
$$
gibt, der obige Sequenz splittet. Da  $\hat{\mathcal{F}}_{(Y,G,(y_i,G_i)^r_{i=0})}$ infinitesimale \"aqui\-va\-riante Deformationen klassifiziert, die auf $Y$ die {\sl triviale} Deformation induzieren, ist sein Tangentialraum genau $\ker(\lambda)$. Deswegen ist 
dieser Tangentialraum isomorph zu $\bigoplus^r_{i=0}D_{\rho'_i}(F[\varepsilon])$, und aus einer Argumentation wie in \cite{BM}, Th.\ 3.3.4 schliesst man, dass $\hat{\Psi}$ formal {\sl \'etale} ist. 

Da obige Funktoren offenbar lokal endlich dargestellt sind, ist auch $\Psi$ als Morphismus von Funktoren  {\'e}tale, d.h.\ die Deformationen, die $${\mathcal{F}}_{(Y,G,(y_i,G_i)^r_{i=0})}$$ klassifiziert, sind lokal-\'etale vollst\"andig bestimmt durch Deformation der Verzweigung in den Verzweigungspunkten.

Die $\hat{\mathcal{F}}_{G_i}$ k\"onnen effektiv pro-dargestellt werden (dies folgt aus \cite{CK2}, Abschnitt 4, insbesondere 4.2.6 und 4.4), und dasselbe gilt f\"ur $\hat{\mathcal{F}}_{(Y,G,(y_i,G_i)_{i=0}^r)}$, da einer
der folgenden F\"alle vorliegt (wie leicht aus der {\sc Riemann-Hurwitz-Zeuthen}schen Formel f\"ur gez\"ugelte Erweiterungen folgt): (i) $g \geq 2$; (ii) $g\geq 1$ und mindestens einen
Punkt ist verzweigt; (iii) $g=0$ und mindestens drei Punkte auf $X$ sind verzweigt ; oder schliesslich (iv) $g=0$ und genau einem Punkt ist wild verzweigt mit Gruppe vom Typ $(1,d)$. 

In den ersten drei F\"allen ist $\Omega_X(D)$ (wobei $D$ der reduzierte Verzweigungsdivisor ist) ampel, und kann man den {\sc Grothendieck}schen Existenzsatz anwenden (\cite{EGA} III). 
Im Falle (iv) gilt dieses nicht, aber da ist $\mathcal{O}_{\P^1}(2)$ ampel. 

\vv

\paragraph \label{linie} {\bf Eine Deformation.} \ Der Einfachheit halber werden wir als Notation f\"ur "`den Funktor dargestellt durch ein Schema $S$"' die Notation "`$S$"' benutzen.  Wir betrachten jetzt den lokalen Funktor $\mathcal{F}_\infty = \mathcal{F}_{G_0}$. Im Abschnitt \ref{familie} haben wir einen Morphismus 
$\til{\mathcal{S}} \rightarrow \mathcal{F}_\infty$ konstruiert. Er ist definiert durch die Familie  $\til{\X} \rightarrow \til{\mathcal{S}}$ aus \ref{subpara-genext4}, weswegen sich durch Komposition ein Morphismus $\Spec \til{\mathfrak{R}}_d \rightarrow \mathcal{F}_\infty$ ergibt. 

Nehmen wir an, das Element $$(\pi: X \rightarrow Y) \in \mathcal{F}_{(Y,G,(y_i,G_i)_{i=0}^r)}(k)$$ (wobei $G_i$ die Verzweigungsgruppen von $\pi$ sind) hat als lokales Bild
$$ \Psi(\pi: X \rightarrow Y) = (\phi_0, \ldots, \phi_r) \in \prod_{i=0}^r \mathcal{F}_{G_i}(k). $$
Wir haben eine  Abbildung 
$$\Spec \til{\mathfrak{R}}_d (k) = \til{\mathbf{L}}_d(k) \rightarrow \mathcal{F}_\infty(k),$$ die surjektiv ist wegen \ref{cor-isomaction1}; sei $\til{\phi}_0$ ein Urbild von $\phi_0$. 

Man w\"ahle wie folgt eine Kurve $$\ell \hookrightarrow \Spec \til{\mathfrak{R}}_d,$$ die durch $\til{\phi}_0$ geht (sei $\infty' \in \ell$ der zugeh\"orige Punkt). Man nehme $u_1 \in V = \theta_{\phi_0}(H) \subset k$ geh\"orend zu $U_1$ in $\til{\mathfrak{R}}_d = k[U_1,\ldots,U_d,W^{-1}]$ und zerlege das zugeh\"orende $\F_q$-lineare Polynom nach Lemma \ref{lem-flagfunctor1} in $$ P^*(S_1,\ldots,S_{d-1}) \circ P^*(S). $$ In dieser  Zerlegung l\"asst man jetzt die  Gr\"o{\ss}en $S_1,\ldots,S_{d-1}$ konstant, w\"ah\-rend $S$ variabel ist (man bemerke nochmals, dass wir $d>1$ voraussetzen). Es ergibt sich dadurch eine Kurve in $\Spec \mathfrak{R}_{(1,d-1)}$ deren Urbild in $\Spec \til{\mathfrak{R}}_d$ wir mit $\ell$ bezeichnen. Man bemerke, dass eine solche Kurve $\ell$ die "`Grenze"' von $\Spec \til{\mathfrak{R}}_d$ in einem allgemeinen Punkt  in der Hyperebene $U_1=0$ ber\"uhrt (wobei "`allgemeiner Punkt"' bedeutet, dass der Punkt nicht zu einer weiteren Grenzkomponente geh\"ort). 

Man konstruiert dann
einen Morphismus $$ \ell \rightarrow \prod_{i=0}^r \mathcal{F}_{G_i}, $$
wobei die Einschr\"ankung $\ell \rightarrow \mathcal{F}_{G_0}$ sich aus der oben konstruierten Kurve ergibt durch 
die Komposition $$\ell \rightarrow \Spec \til{\mathfrak{R}}_d \rightarrow \mathcal{F}_\infty = \mathcal{F}_{G_0},$$ und f\"ur $i>0$ und $S \in \mathbf{Sch}_k$,  $\Hom(S,l) \rightarrow \mathcal{F}_{G_i}(S)$ gegeben wird durch die konstante Abbildung $[ \phi_i \times_k S]$. 

\vv

\paragraph \label{linie2} {\bf Globalisation.} \ Wir behaupten jetzt, dass es eine \'etale Umgebung $S' \rightarrow \ell$ von $\infty'$ gibt, so dass 
$$
\begin{array}{ccc}
S'&\longrightarrow&\mathcal{F}_{(Y,G,(y_i,G_i)_{i=0}^r)}\\
\bigdownarrow & & \bigdownarrow \Psi\\
\ell &\longrightarrow&\prod \mathcal{F}_{G_i}
\end{array}
$$
kommutiert. 
In der Tat, sei $\ell_{/\infty'}$ die formale Komplettierung (als Schema) in $\infty'$. Im formalen Diagramm 
$$
\begin{array}{ccc}
&&\hat{\mathcal{F}}_{(Y,G,(y_i,G_i)_{i=0}^r)}\\
 & & \bigdownarrow \hat{\Psi}\\
\ell_{/\infty'} &\longrightarrow&\prod \hat{\mathcal{F}}_{G_i}
\end{array}
$$
ist $\hat{\Psi}$ formal \'etale (\ref{formal-funktor}).  

Die Eigenschaften der obigen Funktoren aus \ref{formal-funktor} implizieren, dass man das Approximationsverfahren von {\sc M.\ Artin} benutzen kann, um $\ell_{/\infty'}$ durch seine
Henselisierung $\ell_{/^h \infty'}$ zu ersetzen, und es ergibt sich daher der ge\-w\"unsch\-te $S'$ und
\unitlength0.7pt
\begin{center}
\begin{picture}(100,45)
\put(15,35){$\Pi$}
\put(-17,25){$\X'$}
\put(4,30){\vector(1,0){34}}
\put(40,25){$\mathcal{Y}'=S' \times Y$}
\put(0,20){\vector(1,-1){15}}
\put(42,20){\vector(-1,-1){15}}
\put(15,-10){$S'$}
\end{picture}
\end{center}

\vv

\paragraph \label{proof-deg} {\bf Eine Entartung.} \ Das Bild von $S \rightarrow \ell $ ist Zariski dicht und offen, wir k\"onnen also $S'$ um einen Punkt zu $S=S' \cup \{ 0 \}$ erweiteren, derart dass $0$ mit $u_1=0$ korrespondiert (insbesondere wegen der Wahl von $\ell$).  Nachdem, was wir vorher bewiesen haben \"uber $S'=S-\{0\}$, gelten sofort die Eigenschaften (i) und (ii) des Satzes. Wir bilden ein Diagram
 $$
\begin{array}{ccc}
\X'&\longrightarrow&\X''\\
 G \bigdownarrow & & \bigdownarrow \Phi \\
S' \times Y &\longrightarrow&S \times Y
\end{array}
$$
wobei $\X''$ definiert ist als die Normalisierung von $\X'$. Da $\X' \rightarrow S' \times Y$ $G$-Galoissch ist, gilt das gleiche f\"ur $\Phi$. Es ist jetzt aber $\X''$ nicht regul\"ar. In der Tat ist die einzige singul\"are
Faser $\X''_{(0,y_0)}$. W\"ahle $x_0 \in \X''_{(0,y_0)}$, dann gibt es 
wegen \ref{prop-deg3} eine birationale Abbildung 
$$ \Spec \mathcal{O} \dashrightarrow \X''_{/x_0} $$
wobei
$$\mathcal{O}:=k[[x,y,t]]/(P^{\ast}(s_1,\ldots,s_{d-1})(t^{q-1}x^{-q}-x^{-1})^n-y^{-1}) \cap k[[x,y,t]],$$
und $t$ eine Uniformisierende ist in $0$ auf $S$. 
Wir behaupten jetzt, dass diese birationale Abbildung ein Isomorphismus ist. Da sie birational ist, reicht dazu folgendes Lemma aus:

\vv

\paragraph \label{O-ganzab} {\bf Lemma.} \ {\sl $\mathcal{O}$ ist ganzabgeschlossen.} 

\vv

\pf $\mathcal{O}[t^{\pm 1}]$ ist eine $G_{n,d}$-Galoissche Erweiterung von $k[[y,t^{\pm 1}]]$, und deswegen ganzabgeschlossen. Es sei $F$ ganz \"uber $\mathcal{O}$ im Quotientenk\"orper von $\mathcal{O}$, dann ist also
$F = F_0/t^e$ mit $F_0 \in \mathcal{O}$. Da $F$ eine monische Gleichung $F^N + a_{N-1} F^{N-1} + \ldots + a_0 = 0$ erf\"ullt mit $a_i \in \mathcal{O}$, ist $F_0^N \in t \mathcal{O}$. Nun ist $t \mathcal{O}$ ein Primideal, weil $$\mathcal{O}/t\mathcal{O} = k[[x,y]]/(P^{\ast}(s_1,\ldots,s_{d-1})(x^{-1})^n-y^{-1}) \cap k[[x,y]]$$ ein diskreter Bewertungsring ist (\cite{Ser} I Prop.\ 17). Es muss also $F_1 \in \mathcal{O}$ existieren mit $F_0 = t F_1$, und man findet 
$F= F_1/t^{e-1}$ und induktiv $F \in \mathcal{O}$. \qed 

\vv

Die Entartung wird also vollst\"andig von der Proposition \ref{prop-deg3} beschrieben, wobei wir bemerken, dass in der dortigen Beschreibung die entartete Faser in der Tat durch Aufblasen $ x \leadsto t\frac{Z}{W}$ erhalten wird. 

Durch Aufblasen von $S \times Y$ in $(0,y)$ ergibt sich $\mathcal{Y}$, und durch gleichzeitiges Aufblasen von $\mathcal{X}''$ in allen Punkten \"uber $(0,y)$ ergibt sich $\mathcal{X}$.  

Wir zeigen jetzt, dass die horizontalen Komponenten glatt sind. Die Punkte in $\Pi^{-1}(0,y)$ sind glatt wegen \ref{prop-deg3}. Man nehme an, $\xi$ sei singular $\notin \Pi^{-1}(0,y)$. Es wird durch $\Pi$ auf ein Element von $\{0\} \times Y$ abgebildet, aber da $\Pi$ generisch unverzweigt ist, m\"usste $\xi$ dann ein Verzweigungspunkt sein. Weil alle derartigen Verzweigungspunkte durch konstante Deformation der Verzweigungspunkte \"uber $B$ kommen, sind sie insbesondere glatt. 

Die Eigenschaften (iii) und (iv) folgen jetzt sofort.

F\"ur (v) bemerken wir, dass wenn man $\X_{S-\{0\}}$ als gepunktete Kurve betrachtet durch Auszeichnen der Verzweigungspunkte, $\X_{S-\{0\}}$ eine stabile Kurve ist, au{\ss}er
wenn der Ausnahmefall wie in (v) vorliegt. Wenn man obige Konstruktion von $\X_0$ durch Aufblasen  betrachtet, sieht man, dass die gepunktete Struktur
sich auf $\X_0$ \"ubertr\"agt. Dabei haben die vertikalen rationalen Kurven
mindestens drei markierte Punkte, n\"amlich die $q$ Schnittpunkte mit den horizontalen Kurven und zus\"atzlich die Punkte, die von der ganzen Zerlegungsgruppe festgehalten werden ($\{ W/Z=\infty\}$ aus \ref{prop-deg3}). Die horizontalen Kurven bleiben nat\"urlich stabil gepunktet durch die Verzweigungspunkte (au{\ss}er im erw\"ahnten Ausnahmefall); dies ist nur nicht ganz offensichtlich, falls eine horizontale Kurve rational ist. Da begr\"undet man die Stabilit\"at aber folgenderma{\ss}en: da nach \ref{prop-deg3} die Gruppe $G$ die horizontalen Kurven transitiv vertauscht, sind also alle Kurven rational. Da $\X_0$ zusammenh\"angend ist (als flacher Limes von zusammenh\"angenden Kurven), muss eine (ergo alle) solche horizontale(n) Kurve(n) mindestens zwei der vertikalen rationalen Kurven schneiden (die vertikalen Kurven schneiden sich durch die Konstruktion via Aufblasungen nat\"urlich nicht). Deswegen liegen auf jeder der horizontalen Kurven mindestens zwei ausgezeichnete Punkte \"uber $y$ und zus\"atzlich mindenstens ein weiterer ausgezeicheter Punkt, n\"amlich \"uber einem verzweigten Punkt aus $Y-\{y\}$. Deswegen sind auch die horizontalen rationalen Kurven stabil, da mindestens dreifach gepunktet.  \qed

\vv

\paragraph \label{scholie} {\bf Bemerkungen.} \ (i) Im Allgemeinen entartet jede nicht-isotriviale Familie Galoisscher \"Uberlagerungen, die verzweigt \"uber eine {\sl feste} Menge Punkten: siehe {\sc Pries} \cite{Pri2} 3.3.2. Dieses Resultat wurde hier jedoch im Beweis nicht benutzt.

(ii) In \cite{CK2}, Abschnitt 4, wurde bewiesen, dass der infinitesimale Deformationsfunktor einer gez\"ugelten Gruppenoperation vom Typ $(n,d)$ auf $k[[x]]$ pro-dargestellt wird 
durch ein formales Schema vom Typ 
$$ k[[x_0,\ldots,x_{d-1}]]/\langle x_0^{\frac{p-1}{2}} \rangle .$$
Hierbei sind die $d-1$ unbehinderten Deformationsparameter $x_1,\ldots,x_{d-1}$
genau gegeben durch Deformation der Vektorraumeinbettung $ H \hookrightarrow k$. Wenn also bei einer globalen Deformation einer gez\"ugelten Galois\"uberlagerung die Verzweigungspunkte fest bleiben, und die Operation der
Verzweigungsgruppe sich nur in einem Punkt \"andert, so muss sie durch eine solche Deformation von $H \hookrightarrow k$ entstehen. Wir haben in Satz A den Fall, wo ein Unterraum der Dimension eins verschwindet, untersucht. Hieraus leitet man den allgemeinen Fall induktiv ab. 

\vv

\vv

\sectioning{Automorphismen gez\"ugelter Kurven}\label{aut}

\vv

\paragraph \label{intro-aut} {\bf Automorphismen.} \ Wir betrachten jetzt eine Kurve $X$ des Geschlechts $g \geq 2$ \"uber einen algebraisch abgeschlossenen K\"orper der Charakteristik $p>0$. Dann ist $\Aut(X)$ endlich ({\sc F.K.\ Schmidt}), und wir nehmen an, in der Quotientenabbildung 
$$ \pi_X : X \rightarrow Y := \Aut(X) \backslash X$$ sei alle Verzweigung gez\"ugelt (wir nennen dann $X$ gez\"ugelt). {\sc S.\ Nakajima} (\cite{Nak}) hat gezeigt, dass in diesem Falle die {\sc Stichtenoth}sche obere Schranke
$$|\Aut(X)| \leq 16 g^4 + 56 g^3 + 32 g^2 + 4g + 4 \sqrt{1+8 g}(4g^3+4g^2+g) $$
f\"ur die Anzahl der Automorphismen einer beliebigen Kurve $X$ (siehe \cite{Sti}) sich
verbessert zu
$$ |\Aut(X)| \leq 84 g (g-1). $$
(Dies wird dort nur f\"ur ordin\"are Kurven angegeben, aber man \"uberzeugt sich leicht davon, dass der Beweis nur die Tatsache, dass $G_{x,2}=0$ f\"ur alle $x \in X$ benutzt.)

\vv

\paragraph\label{RHZ} {\bf Reduktion.} \ Beide Schranken ergeben sich durch geschickte Minimierung in der {\sc Riemann-Hurwitz-Zeuthen}schen Formel (wo im Falle wilder Verzweigung auch die
h\"oheren Verzweigungsgruppen eine Rolle spielen). Im Laufe der Absch\"atzungen bemerkt man, dass sogar die klassische Hurwitzsche lineare Schranke $|\Aut(X)| \leq 84(g-1)$ gilt, 
au{\ss}er wenn $Y\cong \P^1$, es mindestens einen wildverzweigten Punkt gibt, und die Verzweigung entweder \"uber zwei Punkten liegt, oder $p \neq 2$ und die Verzweigung \"uber drei Punkten
mit Indizes $(2,2,\cdot)$ liegt (\cite{Nak}, p.\ 600). Man kann zus\"atzlich annehmen, dass in wildverzweigten Punkten $x$ die Ungleichung $n_x \neq 1$ gilt, da sonst erneut eine lineare Schranke in $g$ bewiesen werden kann (\cite{Nak}, pp.\ 601--602).

Inwiefern ist das {\sc Nakajima}sche Ergebnis scharf? Wie bereits in der Einf\"uhrung bemerkt, ist keine
Familie $\{X_i\}$, deren Geschlechter $g_i$ strikt steigend sind, bekannt, f\"ur die $|\Aut(X_i)| > \tilde{f}(g_i)$, wobei
$$ \tilde{f}(g) = \max \{ 84(g-1), 2 \sqrt{g} (\sqrt{g}-1)^2 \}. $$
Wir ben\"otigen bald das folgende in diesem Zusammenhang positive Resultat:

\vv

\paragraph \label{Mumford} {\bf Satz} (\cite{CKK}) \ {\sl Sei $K$ ein kompletter nicht-archimedisch
bewerteter K\"or\-per der Charakteristik $p>0$, und sei $X$ eine Mumfordkurve vom Geschlecht $g \geq 2$ \"uber $K$ (d.h.\ die stabile Reduktion von $X$ ist eine Vereinigung von rationalen Kurven, die sich in \"uber den Restklassenk\"orper rationalen Punkten schneiden). Dann gilt 
$|\Aut(X)| \leq \tilde{f}(g)$, und diese obere Schranke ist scharf. }

\vv

Man bemerke, dass Mumfordkurven gez\"ugelt sind im obigen Sinne (Siehe \cite{CKK}, 1.2). 
\vv

\paragraph \label{def-immobil} {\bf Definition.} Wir nennen $X$ {\sl immobil}, 
falls 

(i) $X$ gez\"ugelt ist, 

(ii) $Y =  X / \Aut(X) \cong \P^1$, 

(iii) die Verzweigung entweder \"uber zwei Punkten liegt, oder $p \neq 2$ und die Verzweigung \"uber drei Punkte
mit Indizes $(2,2,\cdot)$ liegt, 

(iv) es mindestens einen wildverzweigten Punkt $x$ gibt mit $n_x \neq 1$,

(v) $s_x = t_x$ in jedem Verzweigungspunkt. 

\vv

\paragraph \label{bew-B} {\bf Beweis von Satz B.} \ Wenn $X$ immobil ist, gilt die obere Schranke wegen der Annahme. Falls eine der ersten vier Annahmen in der Definition (\ref{def-immobil}) nicht erf\"ullt ist, so gilt sogar eine obere Schranke linear im Geschlecht $g$, nach (\ref{RHZ}). Wir k\"onnen also annehmen, es gibt $x$ wildverzweigt mit $s_x \neq t_x$. Wir betrachten eine Familie wie im Satz A. Seien $\{ X_i \}_{i=1}^r$ die gez\"ugelten Kurven in $\X_0 - \Pi^{-1}(E)$. Aus der lokalen Beschreibung \ref{prop-deg3} folgt, dass $\Aut(X)$ die Menge der horizontalen Kurven $\{ X_i \}$ transitiv vertauscht, weswegen sie alle das gleiche Geschlecht haben, das wir mit $\til{g}$ bezeichnen.  Sei $\delta_i$ die Anzahl Knoten von $\X_0$ auf $X_i$. Wir unterscheiden drei F\"alle:

(I) Wenn $\til{g}=0$, so versehen wir den Funktionenk\"orper der Basis $S$ mit der Bewertung, die zum Punkt $0$ geh\"ort, und nehmen die Komplettierung. Dann liegt in einer formalen Umgebung von $\X_0$ in der \"aquivarianten Familie $\X$ eine
glatte Mumfordkurve (mit vollst\"andig zerfallender stabiler Reduktion $\X_0$), und die obere Schranke mit $\til{f}(g) \cdot (g-1)$ ist erf\"ullt nach (\ref{Mumford}). Da wir annehmen, dass $\til{f} \leq f$, folgt das Ergebnis. 

(II) Sei $\til{g} = 1$. Aus der Theorie der Gruppenoperationen ergibt sich die Formel $$|\Aut(X)| = | \Aut(X) \cdot X_1 | \cdot |\mbox{Stab}_{\Aut(X)} X_1 |; $$
hier besteht die Bahn $\Aut(X) \cdot X_1$ genau aus den horizontalen Kurven, und   da $X_1$ triviale Tr\"agheitsgruppe hat, ist $$\mbox{Stab}_{\Aut(X)} X_1 \subseteq \Aut((X_1,\{P_1,\ldots,P_{\delta_1} \})$$ (letztere Gruppe erneut im Sinne der algebraischen Automorphismengruppe von $X_1$, die aber zus\"atzlich die ausgezeichnete Knoten $P_i$ vertauscht). 
 Es ist $\Aut(X_i)$ eine
Erweiterung
$$ 1 \rightarrow T \rightarrow \Aut(X_i) \rightarrow \Aut_0 (X_i) \rightarrow 1, $$
wobei $\Aut_0 (X_i)$ Automorphismengruppe einer elliptischen Kurve (Kurve vom Geschlecht eins mit markiertem Punkt $0$) ist, also $|\Aut_0 (X_i)| \leq 24$ (\cite{Sil} III.10.1), und $T=X_1(k)$ eine (abelsche) Gruppe von Translationen auf der elliptischen Kurve $(X_i,0)$ ist. In unserem Falle soll $$T \cap \Aut((X_1,\{P_1,\ldots,P_{\delta_1} \})$$ sogar die $\delta_1$ ausgezeichnete Punkte vertauschen, weswegen $$|T \cap \Aut((X_1,\{P_1,\ldots,P_{\delta_1} \})| \leq \delta_1.$$ Wir finden also
$$|\Aut(X)| \leq 24 \delta_1 r. $$
Sei $\delta$ die totale Anzahl Knoten auf $\X_0$. Die einzigen Knoten auf den $X_i$ sind die Schnittpunkte mit den vertikalen Kurven, also ist $\delta = r \delta_1.$ 
 Letzten Endes ist 
$$ |\Aut(X)| \leq 24 \delta. $$
Da $\X_0$ gepunktet stabil ist, impliziert die stabile Geschlechterformel, dass 
$ \delta \leq 3g-3, $ (\cite{HM}, p.\ 50) woraus sich dann die lineare obere Schranke $72(g-1)$ ergibt. 

(III) Wenn schliesslich $\til{g} \geq 2$, k\"onnen wir induktiv annehmen, dass $$|\Aut(X_1)| \leq f(\til{g}).$$ Erneut hat $X_1$ triviale Tr\"agheitsgruppe, und also ist 
$$|\Aut(X)| = | \Aut(X) \cdot X_1 | \cdot |\mbox{Stab}_{\Aut(X)} X_1 |. $$
mit $\mbox{Stab}_{\Aut(X)} X_1 \subseteq \Aut(X_1)$. 
Sei $n_i$ die Anzahl der Kurven in $\X_0$ des Geschlechtes $i$. Die nodale 
Geschlechterformel (\cite{HM}, 3.1) impliziert
$$ g-1 = n_{\til{g}}(\til{g}-1) - n_0 + \delta, $$
wobei $\delta$ die Anzahl der Knoten in $\X_0$ ist. Offensichtlich ist $-n_0 + \delta \geq 0$, da jede Kurve mindestens einen Knoten enth\"alt. Deswegen ist
$ g-1 \geq n_{\til{g}} (\til{g}-1)$. Da die Bahn von $X_1$ unter $\Aut(X)$ aus h\"ochstens $n_{\til{g}}$ Kurven besteht, ist $|\Aut(X) \cdot X_1 | \leq \frac{g-1}{\til{g}-1}$ und also 
$$|\Aut(X)| \leq  \frac{g-1}{\til{g}-1} \cdot f(\til{g}). $$
Da wir annehmen, dass $\til{f}(g)/(g-1)$ steigend ist und $g \geq \til{g}$, folgt daraus
schlie{\ss}lich $|\Aut(X)| \leq f(g). $
\qed

\vv

\vv


\setcounter{section}{0}
\renewcommand{\thesection}{\Alph{section}}

\sectioning{\bf Kanonische Form gez\"ugelter Erweiterungen \"uber henselsche Ringe.} \label{hens-norm}

\vv

Wir zeigen jetzt, dass Proposition \ref{prop-canonicalform} auch wahr ist, wenn man statt "`$A$ ist komplett"' nur "`$A$ ist henselsch"' voraussetzt. Wir k\"onnen nach \ref{prop-canonicalform} in den Komplettierungen $\hat{B} = k[[\hat{x}]]/\hat{A}=k[[\hat{y}]]$ eine kanonische Form voraussetzen:
$$ \gamma(\hx) = \zeta \hx, \ \ \ \sigma(\hx) = \frac{\hx}{1-\theta(\sigma) \hx} $$
wobei nach \ref{prop-invariant} $$P^*(\hx^{-1})^n=\hat{y}^{-1}$$
mit $P^*=P^*(a_0,\ldots,a_d)$ f\"ur $a_i \in k$. Sei $k[[\hx]]^H = k[[\hz]]$ mit
$P^*(\hx^{-1})=\hz^{-1}$, $\hz^n=y$ und sei die Operation von $T$ gegeben durch $\gamma(\hz) = \zeta \hz$. 

Man w\"ahle $y \in A$ mit $y = \hy \alpha$ in $\hat{A}$ und $\alpha \in 1 + \mathfrak{m}^2_{\hat{A}}$ (m\"oglich weil $A \hookrightarrow \hat{A}$ dicht ist). Man l\"ose mit Hilfe vom henselschen Lemma ($(n;p)=1$) die Gleichung $\beta^n = \alpha$ mit $\beta$ eindeutig bestimmt durch $\beta \in 1 + \mathfrak{m}^2_{\hat{A}}$.  Sei $z=\hz \beta$, dann ist $z^n=y$ und $\gamma(z)=\zeta z$, weil $\beta \in \hat{A} = \hat{B}^G$. Man definiert dann 
$v \in \mathfrak{m}_{k[[z]]}$ durch $$z=\frac{\hz}{1-\hz v}$$ und l\"ost  die Gleichung $P^*(u)=v$ mit Hilfe vom henselschen Lemma ($a_d \neq 0$) eindeutig nach $u$ in $\mathfrak{m}_{k[[z]]}$. Wenn man $$ x = \frac{\hx}{1-\hx u}$$ definiert, dann ist $x^{-1} = \hx^{-1} - u$, und also 
$$ P^*(x^{-1}) = P^*(\hx^{-1} - u) = P^*(\hx^{-1}) - P^*(u) = \hz^{-1} - v = z^{-1}, $$
weswegen auch $P^*(x^{-1})^n = y^{-1}$. 

Wir berechnen jetzt $\gamma(x)$. Weil $$\gamma(z)=\zeta z = \zeta \frac{\hz}{1-\hz v}$$ gleich $$\frac{\gamma(\hz)}{1-\gamma(\hz) \gamma(v)} = \frac{\zeta \hz}{1-\zeta \hz \gamma(v)}$$ ist, findet man $\gamma(v)=\zeta^{-1} v$. Weil $\gamma(P^*(u)) = \gamma(v) = \zeta^{-1} v$ gleich $P^*(\gamma u)$ ist, muss $$P^*(u) = P^*(\zeta \gamma(u))$$  gelten, und weil $u$ eindeutig ist in $\mathfrak{m}_{k[[z]]}$ und $\zeta \gamma(u) \in \mathfrak{m}_{k[[z]]}$ (weil $G$ als lokaler Morphismus auf $B$ operiert und $\zeta \in \F^\times_q$), muss dann $\gamma(u) = \zeta^{-1} u$ sein. 
Daraus ergibt sich, dass $$\gamma(x) = \frac{ \gamma(\hx)}{1-\gamma(\hx) \gamma(u)} = \frac{\zeta \hx}{1-\zeta \hx \zeta^{-1} u} = \frac{\zeta \hx}{1-\hx u} = \zeta x.$$ 

Was $\sigma(x)$ angeht, so haben wir $$\sigma(x) = \frac{ \sigma(\hx) }{1-\sigma(\hx) \sigma(u)} $$ mit $u \in k[[z]] = \hat{B}^H$, also $\sigma(u)=u$. Mit Hilfe der obigen expliziten Formeln berechnet man sofort, dass 
$$\sigma(x) = \frac{x}{1-\theta(\sigma) x}.$$

Wir haben $y \in A$ gew\"ahlt, und $x \in \hat{B}$ erf\"ullt $P^*(x^{-1})^n=y^{-1}$. Weil $A$ henselsch ist (und deswegen auch $B$), kann man $x' \in B$ finden mit
$P^*(x'^{-1})^n =y^{-1}$ und $x \equiv x' \mbox{ mod } \mathfrak{m}_B^N$, und da die Gleichung nur endlich viele L\"osungen hat, findet man f\"ur $N \rightarrow \infty$, dass $x \in B$. \qed

\vv

Man \"uberzeugt sich jetzt leicht davon, dass auch  \ref{prop-isomaction}, \ref{prop-invariant} und \ref{cor-isomaction1} war sind, falls man nur "`$A$ henselsch"' voraussetzt.

\vv

\vv

\sectioning {\bf Konfigurierende Funktoren} \label{anhang}

\vv

In \cite{Pri2}, 2.2.1 hat {\sc Pries} einen Begriff "`Konfigurationsraum"' eingef\"uhrt, den wir jetzt leicht erweitern.

\vv

\paragraph \label{def-konfig} {\bf Definition.} Es seien $\mathcal{F}$ und $\mathcal{G}$ zwei kontravariante Funktoren $\mathbf{Sch}_k \rightarrow \mathbf{Mengen}$. Wir sagen, dass der Funktor $\mathcal{G}$ den Funktor $\mathcal{F}$ {\sl konfiguriert}, falls es einen Morphismus $$ \mathcal{F} \stackrel{\Psi}{\longrightarrow} \mathcal{G}$$ 
 zwischen Funktoren gibt, wobei gilt:

(i) F\"ur jeden algebraisch abgeschlossenen K\"orper $F/k$ ist $\Psi(F)$ bijektiv;

(ii) F\"ur jedes Schema $S \in \mathbf{Sch}_k$ und jedes $\phi \in \mathcal{F}(S)$ gibt es $S' \in \mathbf{Sch}_k$, eine endliche radiziale Abbildung $i \colon S' \rightarrow S$ und ein eindeutig bestimmtes $\phi' \in \mathcal{G}(S')$ mit $\Psi(\phi') = i^*(\phi)$. 

\vv

\paragraph \label{bem-konfig} {\bf Bemerkung.} \ Falls der Funktor $\mathcal{G}$ den Funktor $\mathcal{F}$ konfiguriert und $\mathcal{G}$ fein dargestellt wird durch ein Schema $T$, dann ist $T$ ein Konfigurationsraum f\"ur $\mathcal{F}$ im Sinne von \cite{Pri2}, 2.2.1. 

\vv

\paragraph \label{prop-konfig} {\bf Proposition.} {\sl Sei $A$ ein henselscher diskreter Bewertungsring \"uber einen algebraisch abgeschlossenen K\"orper $k$ mit Restklassenk\"orper $k$, sei $G=G_{n,d}$, und sei $\mathcal{F}_{A,G}$ der Kontravariante Funktor $\mathbf{Sch}_k \rightarrow \mathbf{Mengen}$ definiert
durch 
$$S\mapsto\mathcal{F}_{A,G}(S)=
\left\{
\begin{minipage}{17em}
\setlength{\baselineskip}{.85\baselineskip}
\begin{small} \begin{center}
{\slshape $G$-\"Uberlagerung $\X \rightarrow \Spec A \times_k S$, die in jedem geometrischen Punkt gez\"ugelt sind}
\end{center} \end{small}
\end{minipage}
\right\} / \cong
$$
wobei geometrischer Punkt separabeler Abschluss eines Punktes bedeutet. Dann konfiguriert $\mathcal{L}_d/\mathbf{G}_{m,k}$ den Funktor $\mathcal{F}_{A,G}$.}

\vv

\pf Wir konstruieren zuerst einen Morphismus $\mathcal{L}_d / \mathbf{G}_{m,k}  \stackrel{\Psi}{\longrightarrow} \mathcal{F}_{A,G}$. Dazu bemerken wir, dass im affinen Fall $\til{L}_d/\mathbf{G}_{m,k} \cong \til{\mathcal{S}}$ die universelle Familie $\til{\mathcal{X}}$ mit $G$-Operation tr\"agt (\ref{subpara-genext3}, \ref{subpara-genext4}). Man bekommt also offensichtlich einen Morphismus
$$ \til{\mathcal{S}} \rightarrow \mathcal{F}_{A,G} |_{\mathbf{Aff}^o_k}, $$
(wobei ${\mathbf{Aff}^o_k}$ die Kategorie der affinen Schemata \"uber $k$ ist),  der f\"ur die triviale Operation von $\mathfrak{G}$ auf der rechten Seite $\mathfrak{G}$-\"aquivariant ist. Deswegen ergibt sich durch Garbifizierung 
ein Morphismus $$\mathcal{L}_d / \mathbf{G}_{m,k} = \mathfrak{G} \backslash \til{\mathcal{L}}_d / \mathbf{G}_{m.k} \rightarrow \mathcal{F}_{A,G}.$$ 

Eigenschaft (i) ist klar wegen \ref{cor-isomaction1}. F\"ur Eigenschaft (ii) w\"ahlen wir $S \in \mathbf{Sch}_k$ und setzen $S':=S_{\mathrm{red}}$ (das reduzierte Schema von $S$). Offenbar ist $S' \rightarrow S$ endlich radizial, und wir zeigen jetzt, dass es einen Isomorphismus gibt
$$ (\mathcal{L}_d / \mathbf{G}_{m,k})(S') \stackrel{\sim}{\longrightarrow} \mathcal{F}_{A,G}(S').$$ 
Wir k\"onnen dabei annehmen, dass $S'$ reduziert ist und sogar affin: $S'=\Spec R$. 

Die Abbildung ist surjektiv: Sei $$\phi=(\mathcal{X} \stackrel{G}{\longrightarrow} \Spec A \times_k \Spec R) \in \mathcal{F}_{A,G}(S').$$
Durch Komplettierung bekommen wir
 $$
\begin{array}{ccc}
\hat{\mathcal{X}} & {\stackrel{G}{\longrightarrow}} & \Spec R[[y]]\\
\bigdownarrow & & \bigdownarrow  \\
\mathcal{X} & {\stackrel{G}{\longrightarrow}} & \Spec A \times_k \Spec R
 \end{array}
$$
mit $\hat{A} = k[[y]]$. 
Wir schreiben $\bar{s}$ f\"ur den separabelen Abschluss eines Punktes $s \in S'$ und erweitern die obige Gruppenoperation durch Basiswechsel via $\Spec \hat{\mathcal{O}}_{S',\bar{s}} \rightarrow \Spec R$ zu
 $$
\begin{array}{ccc}
\hat{\mathcal{X}} \times_{\Spec R} \Spec \hat{\mathcal{O}}_{S',\bar{s}} & \stackrel{G}{\longrightarrow} & \Spec \hat{\mathcal{O}}_{S',\bar{s}}[[y]] \\ 
 \bigdownarrow & & \bigdownarrow  \\
\hat{\X} & \stackrel{G}{\longrightarrow} & \Spec R[[y]] \end{array}
$$
Wir haben jetzt einen Isomorphismus $$\hat{\mathcal{X}} \times_{\Spec R} \Spec \hat{\mathcal{O}}_{S',\bar{s}} \stackrel{\sim}{\longrightarrow} \Spec \hat{\mathcal{O}}_{S',\bar{s}}[[x]].$$ 
Dies folgt aus (\cite{CK2}, Abschnitt 4), wo der verselle formelle infinitesimale Deformationsring 
berechnet wurde (und man Effektivit\"at hat wie im Beweis des Hauptsatzes A (\ref{linie})). Mit Hilfe der Artinschen Approximation k\"onnen wir $S$ durch eine \'etale Umgebung ersetzen; danach sind wir in der Situation 
$$ \hat{\mathcal{X}} = \Spec R[[x]] \stackrel{G}{\longrightarrow} \Spec R[[y]].$$  Ein Element $\sigma \in H \subseteq G$ operiert auf $R[[x]]$ durch 
$$ \sigma(x) = a_0 +a_1 x + a_2 x^2 + \ldots, $$
und wir haben $a_0=0$ und $a_1=1$ modulo alle 
Primideale $\mathfrak{p}$ von $R$. Da $R$ aber reduziert ist, ist $\bigcap \mathfrak{p} = \{ 0 \}$, weswegen $a_0=0, a_1=1$. Wir k\"onnen also
$\sigma(x) = x + \theta(\sigma) x^2 + \ldots$ schreiben, und $\theta: H \rightarrow R$ bildet $H$ auf einen $\F_q$-Vektorraum in $R$ ab: Durch direkte Berechnung stellt sich heraus, dass $\theta$ additiv ist, und wenn wir 
die Operation $\tau(x) = \zeta x + \ldots$ eines Torusteils ber\"ucksichtigen, finden wir leicht, dass $\zeta^{-1} \theta(H) = \theta(H)$. Desweiteren ist $\theta(H)$ nicht-entartet, da $\theta$ f\"ur jede Spezialisierung injektiv ist. Also bekommen wir ein Element $\theta(H)$ in  $\til{\mathcal{L}}_d$, das auf $\phi$ abgebildet wird. 

Die Abbildung ist auch injektiv: Wenn $V$ und $V'$ in $R$ die gleiche
\"Uberlagerung \"uber $\Spec R \times \Spec A$ liefern, sind die $\omega$-Invarianten gleich in jedem geometrischen Punkt $\bar{s}$ von $S'$ (\ref{prop-isomaction} gilt f\"ur henselsche Ringe), weswegen $V$ und $V'$ sich h\"ochstens um ein Skalar in $R^\times$ unterscheiden. \qed

\vv

In \ref{subpara-genext3} wurde gezeigt, dass das Schema $\til{\mathcal{S}}$ den Funktor
$\til{\mathcal{L}}_d / \mathbf{G}_{m,k}$ grob darstellt. Wir haben gerade gesehen, dass dieser Funktor $\mathcal{F}_{A,G}$ konfiguriert. Man kann also $\til{\mathcal{S}}$ als eine Art von "`grobem Konfigurationsraum"' f\"ur $\mathcal{F}_{A,G}$ betrachten. 

\vv

\vv

\begin{small}

\def\refname{\noindent\normalsize{Literatur}}

\vv

\vv

\noindent Fachbereich Mathematik,
Universit\"at Utrecht, Postfach 80010, 3508 TA  Utrecht, Niederlande (gc)

\vv

\noindent Universit\"at Kioto, Fakult\"at der Wissenschaften, Fachbereich Mathematik, Kioto 606-8502, Japan (fk)

\vv

\noindent Email: {\tt cornelissen@math.uu.nl}, {\tt kato@math.kyoto-u.ac.jp}

\end{small}

\end{document}